\documentclass [12pt,a4paper,reqno]{amsart}
  \input amssymb.sty

\usepackage{epsfig}
\usepackage{youngtab}

\newcommand{\ef}{\end{equation}}
\chardef\bslash=`\\ 

\newtheorem{thm}{Theorem}[section]
\newtheorem*{thm*}{Theorem}

\newtheorem{lem}{Lemma}[section]
\newtheorem{corl}{Corollary}[lem]

\newtheorem{prop}{Proposition}[section]
\newtheorem{prop*}{Proposition}

\theoremstyle{definition}
\newtheorem{defn}{Definition}[section]
\newtheorem{examp}{Example}
\newtheorem*{examp*}{Example}
\newtheorem*{remark*}{Remark}
\newtheorem*{CC*}{Crossover Conjecture}
\newtheorem*{Note*}{Note}
\newtheorem*{defn*}{Definition}

 \theoremstyle{remark}
\newtheorem{remark}{Remark}[section]

\newcommand{\wt}{\widetilde}

 \renewcommand{\sectionmark}[1]{}

\newcommand{\Irr} {\operatorname{Irr}}

\renewcommand{\a}{\alpha}

 \date{}
\begin{document}

\title  [Crossover Morita Equivalences] {Crossover Morita equivalences of spin blocks of the symmetric and alternating groups}
\author[R. Leabovich, M. Schaps]
{ Ruthi Leabovich  and Mary Schaps}
 \address{Department of Mathematics, Bar-Ilan University, 52900 Ramat-Gan, Israel}
 \email {Ruthi.Leabovich@gmail.com}

 \address{Department of Mathematics\\
 Bar-Ilan University, 52900 Ramat Gan, Israel}
\email{mschaps@macs.biu.ac.il}

 \begin{abstract} We demonstrate source algebra equivalences between spin blocks of the families of covering groups $\{\tilde S_n\}$ and $\{\tilde A_n\}$ of the symmetric and alternating groups, for pairs of blocks at the ends of maximal strings.  These equivalences remain within the family of groups if the cores of the two blocks have the same parity and cross over from one family to the other if the cores are of opposite parity.  This demonstrates the Crossover Conjecture of Kessar-Schaps for the easier case of extremal points of maximal strings.  The problem of  demonstrating crossover derived equivalence between symmetrically-placed blocks in the interior of maximal strings remains open.

 The Brou\'e conjecture, that a block with abelian defect
 group is derived equivalent to its Brauer correspondent, has been proven for
 blocks of cyclic defect group and verified for many other blocks.  This paper is part of a study of the spin block case, but has wider application, since it establishes Morita equivalences also when the defect group is not abelian.  The results of the paper allow us to give a sharp upper bound both for the maximal N required to get representatives of all Morita equivalence classes and also for the number of Morita equivalence classes of spin blocks for a given $p$.  This improves the bounds given in the proof of the Donovan conjecture for spin blocks given by Kessar.

     Finally, we relate the Scopes involutions used in this paper, slightly reworked from those in \cite{K}, with reflections, in the simple roots, of the long roots of a twisted affine Lie algebra.
 \end{abstract}
\thanks{Work is partially supported by a grant from the Bar-Ilan University Research Authority, and formed part of the first author's Ph.D. thesis at Bar-Ilan}

\maketitle

\large
\section{Introduction}

Let $G$ be a finite group and $\mathcal F$ be an algebraically closed field of characteristic $p > 0$. We are primarily interested here in the modular case when $p>0$. Writing
\[ kG=\bigoplus^{r}_{j=1}B_j.\] as a decomposition into blocks, we let $D_i$ be the \textbf{defect group} of the block $B_i$, determined up to conjugacy in $G$. The defect group are the $p$-groups of order $p^d$, and the integer $d$ is called the \textbf{defect} of the block. By a fundamental result of Brauer \cite{A}, there is a one-to-one correspondence between blocks of $\mathcal FG$ with defect group $D$ and blocks of the group algebra of the normalizer $N_G(D)$ with the same defect group $D$. Not only is $N_G(D)$ generally much smaller then $G$, its blocks have a relatively simple structure, because the defect group is normal. The block of $\mathcal FN_G(D)$ corresponding to a block $B$ of $kG$ is called its \textbf{Brauer correspondent}.

Although there are infinitely many blocks with given defect group $D$, Donovan has conjectured that there are only a finite number of Morita equivalence classes.  Puig has generalized this to a conjecture that for a given defect group, there are only finitely many classes of blocks up to Puig equivalence, in which blocks are equivalent if they have the same source algebra.

One can still hope for some weaker sort of equivalence, such as derived equivalence, between blocks. In particular, two blocks are \textbf{derived equivalent} if there is an equivalence of categories between the corresponding \textbf{bounded derived categories}, where the bounded derived category $D^b(B)$ of a block $B$ is the category of bounded complexes of finite dimensional $B$-modules in which quasi-isomorphisms have been formally inverted. If two blocks are derived equivalent, then they share many important invariants, including the number of simple modules. Moreover, they have isomorphic center, Hochschild cohomology, and, to some extent, the same deformation theory.  Brou\'e has conjectured \cite{B1}, \cite{B2}, that if the defect group $D$ is abelian, the block will be derived equivalent to its Brauer correspondent.  This would reduce the number of possible derived equivalence classes to the number of possible blocks with normal defect group $D$, which is generally much smaller than the number of Morita equivalence classes.

We are considering the possible Morita equivalences among blocks in the two families,
$\{\tilde S_n\}$ and $\{\tilde A_n\}$, and we demonstrate that for this purpose the two families should be treated together. The blocks are determined by a core $\rho$ and a weight $w$.  The Crossover Conjecture (Kessar-Schaps) \cite{AS} asserts that for any $w$ there are exactly two derived equivalence classes in the union of the blocks from the two families.

In $\S 2$ we review the definition of the block-reduced crystal graph from \cite{AS}.  In $\S 3$, we give the precise combinatorial criteria for a pair of blocks to lie at the extremal points of a maximal $i$-string in the block-reduced crystal graph.  In $\S4$, we demonstrate that such blocks, with crossovers where necessary, are source algebra equivalent, and thus Morita equivalent.

For any core $\rho$, let $B_{\rho^w}$ be the block algebra of $\mathcal O\widetilde S_n$ with core $\rho$ and weight $w$ and let $B'_{\rho^w}$ be the corresponding block of $\mathcal O\widetilde A_n$ with core $\rho$ and weight $w$.
We will prove the following:
\begin{thm*} Suppose that the blocks with cores $\nu$ and $\mu$ and weight $w$ lie at the ends of a maximal strings in the block reduced crystal graph.  Then
if the parities are the same, $B_{\nu^w}$ is source algebra equivalent to $B_{\mu^w}$, and $B'_{\nu^w}$ is source algebra equivalent to $B'_{\mu^w}$. If the parities are  different,$B_{\nu^w}$ is source algebra equivalent to $B'_{\mu^w}$, and $B'_{\nu^w}$ is source algebra equivalent to $B_{\mu^w}$.
\end{thm*}

  In $\S 5$, we give the sharp bound for Donovan's conjecture and exhibit the block for which the bound is attained. Finally, in $\S 6$, we describe the connection between the Scopes involutions used in $\S3 - \S5$ and the reflections through simple roots of the long roots of a twisted affine Lie algebra of type $A$, using the corresponding parametrization of the blocks to give a bound on the number of possible Morita equivalence classes for a given weight $w$.

\section{Definitions and notation}
\subsection{Spin block of the symmetric and alternating groups.}
The symmetric group $S_{n}$ and the alternating group $A_{n}$ have central extension $\widetilde{S}_{n}$, $\widetilde{A}_{n}$
with kernel $C_{2}$ the cyclic group of order 2. The group algebra of one of the covering groups can be decomposed into two subalgebra of equal dimension by the value of the characters on the central involution $z$. One is isomorphic to the group algebra of the original group.  The characters in the second component, for which $z$ takes the value $-1$, will be called spin representations, and the corresponding blocks will be called \textsl{spin blocks}.

For every block of $S_{n}$ the isomorphism classes of irreducible representations are labeled by the partitions:  $(n),(n-1,1),(n-2,2),(n-2,1,1),... $ For example: in $S_{6}$ irreducible representations are labeled by the partitions  \[(6), (5,1), (4,2), (4,1,1), (3,3), (3,2,1), (3,1,1,1),\]\[ (2,2,2), (2,2,1,1), (2,1,1,1,1), (1,1,1,1,1,1).\]
For every spin block of $\widetilde{S}_{n}$ and $\widetilde{A}_{n}$, the isomorphism classes of the irreducible representations are labeled by the \textsl{strict} partitions $\lambda=(\lambda_{1},...,\lambda_{r})$ where $\lambda_{i}\neq \lambda_{j}$ for $i, j$ satisfying $1\leq i,j \leq r$. For example: in $\widetilde{S}_{6}$, irreducible spin representations are labeled by the strict partitions \[(6),(5,1),(4,2),(3,2,1).\]

\begin{defn} A \textsl{$p$-bar core} is a strict partition that does not contain any parts divisible by $p$, and does not contain parts congruent to $i$ and to $p-i$ for any $i$ satisfying $1\leq i \leq p-1$.
\end{defn}

\begin{defn} \textsl{Removing a $p$-bar}:

The strict partition can be represented by an abacus with $p$ runners labeled by the residues ${0,1,...,(p-1)}$ where the parts of the strict partition are represented as beads. The part $\lambda_{i}=ap+b$ corresponds to a bead in the runner $b$ with height $a$. Removing a $p$-bar consists either of
 \begin{enumerate}
 \item Lowering the position of one bead one place down on its runner into an empty place, which corresponds to reducing a single part $\lambda_{i}$ by $p$ (if $\lambda_{i}-p$ is not a part of $\lambda$),
 \item  Removing the bottom bead in the 0-runner, which corresponds to removing the part $p$ from the strict partition  $\lambda$, or
 \item Removing two beads which sum to $p$.

 In each case we reduce the sum of the parts of the strict partition by $p$.
 \end{enumerate}
\end{defn}
Every row in the character table of  $\widetilde{S}_{n}$  or $\widetilde{A}_{n}$ corresponding
an irreducible spin representation is labeled by a strict partition.
(It can be the same strict partition for a pair of  rows). When one removes the maximal number $w$ of $p$-bars  from the strict partition one arrives at a $p$-strict partition $\rho$ called the $p$-core, and the integer $w$ is called the weight of the block. The $p$-bar core of
a partition $\lambda$ will be represented by $\rho (\lambda)$.
All the representations that have the same $p$-core belong to the same spin block \cite{BK}.

\begin{defn}\cite{AS}
Let $t=\frac{p-1}{2}$.  Any $p$-core $\rho$ can be represented by a \textsl{core t-tuple} \[c(\rho)=((l_{1},\epsilon_{1}),...,(l_{t},\epsilon_{t}))\] where $l_{i}$ is the number of beads on the runner numbered $i$ or $p-i$, and we set $\epsilon_{i}=0$ if there are beads on runner $i$, and otherwise we set $\epsilon_{i}=1$. Note the rather counter-intuitive choice that if there are no beads on either runner, so that $l_i=0$, then $\epsilon_i = 1$. We will abbreviate $c(\rho(\lambda))$ by $c(\lambda)$.
In what follows, the cores $\rho$ will generally be represented by their core $t$-tuple $c(\rho)$, since the description of the actual partition is too bulky and does not exhibit the special properties of the core.

\end{defn}

\newtheorem{def*}{Definition}
\begin{defn} A \textsl{p-strict partition} is a partition $\lambda=(\lambda_{1},...,\lambda_{r})$ is a partition that has no repeated parts except possibly parts divisible by $p$.  A \textsl{p-restricted partition} $\lambda$ is a $p$-strict partition that satisfies the conditions:
\begin{enumerate}
	\item $\lambda_{i}-\lambda_{j}\leq p$,
	\item $\lambda_{i}-\lambda_{j}=p$,  implies that $p$ does not divide $\lambda_{i}$.
\end{enumerate}
The set of $p$-restricted partitions will be denoted by $RP_{p}$. The set of strict partitions will be denoted by $DP$ (for distinct parts). Note that not all elements in $RP_{p}$ are in $DP$, since parts divisible by $p$ can be repeated.
\end{defn}

\subsection{Crystal graphs, the symmetric case}
The algebra $\mathcal F \widetilde{S_{n}}$ contains a commutative subalgebra $\mathcal A_{n}$ generated by the preimages in $\mathcal F \widetilde{S_{n}}$ of the squares of the Jucys-Murphy elements in $k{S_{n}}$\\
$L_{1}=0$\\
$L_{2}=(1 2)$\\
$L_{3}=(1 3)+(2 3)$\\
$L_{4}=(1 4)+(2 4)+(3 4)$\\
\vdots\\
$L_{n}=\sum^{n-1}_{i=1}(in)$\\

The \textsl{Young diagram} of a partition $\lambda=(\lambda_{1},...,\lambda_{r})$ is a set of $r$ rows of boxes, with $\lambda_{i}$ boxes in row $i$. Each row will be filled in, as far as possible, by repetitions of the sequence $0,1,...,t-1,t,t-1,...,1,0$.

For each residue $i$, we let $q(i)=i(i+1)$.
The $q(i)$ are the possible eigenvalues of the preimages of elements $[L^{2}_{j}]$.
Let the $p$-bar residue $\hat{m}$ of an integer $m$ be defined by
\[
\hat{m}=\left\{  \begin{array}{ll}
         \bar{m} & \mbox{$0\leq\bar{m}\leq \frac{p-1}{2}$};\\
         p-1-\bar{m} & \mbox{$\frac{p+1}{2}\leq\bar{m}\leq p-1$}.\end{array} \right.
\]
where $\bar{m}$ is the usual residue mod $p$, $0\leq\bar{m}\leq p-1$. Note that if $0\leq \bar{m}\leq \frac{p-1}{2}$, then $q(\bar{m})=q(p-\bar{m}-1)$. For example if $p=5$ the $p$-bar residues of the integers 0,1,2,3,4,... are 0,1,2,1,0,0,1,2,1,0,...
The \textsl{content} $\gamma(\lambda)$ is the $(t+1)$-tuple $\gamma=(\gamma_{0},...,\gamma_{t})$, where $\gamma_{i}$ is the number of boxes containing the integer $i$.  The content will determine the block to which the irreducible module belongs.

We now give a description of the crystal graph of $\widetilde{S_{n}}$. Let \\$i\in I=\left\{1,2,...,t\right\}$ be some fixed residue. Every node $A \in \lambda$ can be written $A=(r,s)$ where $r$ is number of the row of  the node $A$, and $s$ is the number of the box on the row $r$, therefore $s\leq \lambda_{r}$ and this box is filled with $\hat{s}$.
A node $A=(r,s) \in \lambda$ is called \textsl{i-removable} (for $\lambda$) if one of the following holds \cite{BK}:
\begin{enumerate}
	\item $res A=i$ and $\lambda-{A}$ is also a $p$-strict partition.
	\item the node $B=(r,s+1)$ immediately to the right of $A$ belongs to $\lambda$, $res A=res B=0$, and both $\lambda-\left\{B\right\}$ and $\lambda-\left\{A,B\right\}$ are $p$-strict partitions.
\end{enumerate}
Similarly, a node $B=(r,s)\notin \lambda $ is called \textsl{i-addable} (for $\lambda$) if one of the following holds:
\begin{enumerate}
  \item $res B=i$ and $\lambda \cup \left\{B\right\}$ is again a $p$-strict partition.
  \item the node $A=(r,s-1)$ immediately to the left of $B$ does not belong to $\lambda$, $res A=res B=0$ and both $\lambda \cup \left\{A\right\}$ and $\lambda \cup \left\{A,B\right\}$ are $p$-strict partitions.
\end{enumerate}

Now label all $i$-addable nodes of the Young diagram of $\lambda$ by + and all $i$-removable nodes by -. The $i$-signature of $\lambda$ is the sequence of pluses and minuses obtained by going along the rim of the Young diagram from bottom left to top right and reading of all the signs. The \textsl{reduced $i$-signature} of $\lambda$ is obtained from the $i$-signature by successively erasing all neighboring pairs of the form +-, so that we get a sequence of -'s followed by +'s. Nodes corresponding to a sign - in the reduced $i$-signature are called \textsl{$i$-normal}, while nodes corresponding to a sign + are called \textsl{$i$-conormal}. The rightmost $i$-normal node (corresponding to the rightmost - in the reduced $i$-signature) is called \textsl{$i$-good}, and the leftmost $i$-conormal node (corresponding to the leftmost + in the reduced $i$-signature) is called \textsl{$i$-cogood}.

We now define the the edges in the crystal graph for $\widetilde{S_{n}}$ and $\widetilde{A_{n}}$ by the connecting two $p$-restricted partitions $\lambda$, $\chi$ by an edge labeled $i$ if one is obtained from the other by removing an $i$-good node or adding an $i$-cogood node. If $\lambda$ and $\lambda'$ have weights $w(\lambda) \leq w(\lambda')$, then all of the $p$-restricted partitions in the block of $\lambda$ will map one-to-one into the $p$-restricted partitions in the block of $\lambda'$. Let $e_{i}$ denote this map of partitions in the rank-reducing direction, and $f_{i}$ the map in the rank increasing direction. Both $e_{i}$ and $f_{i}$ are nilpotent, that, a large power is always zero.

\begin{defn}\cite{AS} The \textsl{block-reduced crystal graph} has as vertices the labels $\rho^w$ of blocks, with edges wherever there are edges between some of the simples in the regular crystal graph.  A \textsl{maximal $i$-string} is a maximal sequence of blocks connected by by restriction and induction operators $e_i$ and $f_i$.
\end{defn}

For $0<i \leq t$ we define a \textsl{Scopes involution}
\begin{defn}
Let $D$ be the set of all partitions which are strict, except that parts divisible by $p$ may occur any number of times.  We define involutions of $D$ given by:

$$K_{i}:D \rightarrow D, ~\textsl{for}~i = 0,\dots,t.$$

\begin{itemize}

	\item for $0<i<t$: the involution $K_{i}$  interchanges the beads on runner $i$ and $i+1$ as well as the beads on runner $p-i$ and $p-i-1$. For cores we get that \[c(\lambda)=((l_{1},\epsilon_{1}),...,(\ell_{i},\epsilon_{i}),(\ell_{i+1},\epsilon_{i+1}),...,(\ell_{t},\epsilon_{t}))\] after acting by $K_{i}$, becomes \[c(\lambda')=((\ell_{1},\epsilon_{1}),...,(\ell_{i+1},\epsilon_{i+1}),(\ell_{i},\epsilon_{i}),...,(\ell_{t},\epsilon_{t})).\]
	\item for $i=t$: the involution $K_{t}$  interchanges beads on runner $t+1$ and runner $t$, (recall $t=\frac{p-1}{2}$) and parts divisible by $p$ are remain fixed so on cores we get that \[c(\lambda)=((\ell_{1},\epsilon_{1}),...,(\ell_{t},\epsilon_{t})),\] after acting by $K_{t}$, becomes \[c(\lambda')=((\ell_{1},\epsilon_{1}),...,(\ell_{t},1-\epsilon_{t}))\].
  \item For $i=0$: the involution $K_{0}$ interchanges the beads $ap+1$ above the place for $1$ on runner $1$ with the beads of form $bp-1$ on runner $p-1$, and either adds the part $1$ if is not there or removes it, if it is.
The involution $K_0$ interchanges $D^+$ and $D^-$, where

$D^{+}$=\{  partitions in $D$ with the part $1$\},

  $D^{-}$=\{  partitions in $D$ without the  part $1$ \}.

For cores we get that \[c(\lambda)=((\ell_{1},\epsilon_{1}),\dots,(\ell_{t},\epsilon_{t}))\] after acting by $K_{i}$, becomes \[c(\lambda)=((\ell_{1}-(-1)^{\epsilon_1},1-\epsilon_{1}),\dots,(\ell_{t},\epsilon_{t})).\]

\end{itemize}
\end{defn}

\begin{remark}\label{Scopes} For $i>0$, the involution given here, restricted to strict permutations, is exactly the involution $\widetilde {Sc}_{i+1}$ defined in \cite{K}. For $i=0$, the map $K_0$ is actually an involution, unlike $\widetilde {Sc}_1$ of \cite{K}.
\end{remark}

\begin{remark} In $\S6$, we will demonstrate that the Scopes involution $K_i$ correspond to a reflection in a simple root $\a_i$ of the long roots of an affine Lie algebra.
\end{remark}

\begin{examp}
For $p=5$, let $$\rho=(12,7,6,2,1),\lambda=(12,11,7,6,4,2,1)$$
 be, respectively, the $p$-core and a partition in the block $\rho^3$.  Then
$$K_0(\rho) = (12,7,4,2), K_0(\lambda) = (12,9,7,6,4,2);$$
$$K_1(\rho) = (11,7,6,2,1),K_1(\lambda) =(12,11,7,6,3,2,1);$$
$$K_2(\rho) = (13,8,6,3,1),K_2(\lambda)=(13,11,8,6,4,3,1).$$

\end{examp}

In order for these involutions to be of use, we must show that they preserve blocks, which is equivalent to showing that they preserve cores.

\begin{lem} If $\lambda \in D$, and if $K_i(\lambda) = \chi$, then $$K_i(\rho(\lambda)) = \rho(\chi).$$
\end{lem}
\begin{proof}

 Case 1: $\lambda \in DP$, i.e., all the parts of $\lambda$ are distinct.  For $i > 0$, this was done in \cite{K}, Lemma 4.7, with the changes in notation given in Remark \ref{Scopes}.  It remains to demonstrate it for $K_0$.  If $\lambda \in D$, with weight $w$, let $\rho = \rho(\lambda)$ be its core, and let $\chi = K_0(\lambda)$. We must show that $\rho'=\rho(\chi)$ is equal to $K_0(\rho)$. Obviously, since all the beads on runners other than $1$ and $p-1$ are identical, we can make the same moves both for $\lambda$ and for $\chi$, with identical results, and since the order of the moves is irrelevant to the final result, we assume that these moves have been made first.  Therefore, we need only analyze the moves on runners $1$ and $p-1$.

 Let $A=\{a>0|ap+1 \in \lambda \}$ and $B=\{b>0|bp-1 \in \lambda \}$.  The effect of performing $K_0$ is to interchange $A$ and $B$, while also adding or removing the part $1$. The first thing we do, in either $\lambda$ or $\chi$, is to close up $A$ and $B$ into intervals $A_0 = \{1,\dots,|A|\}$ and $B_0=\{1,\dots,|B|\}$, a procedure which does not effect either core and which takes the same number of moves for either partition.

 \textbf{Subcase 1.a:} \textsl{$\lambda \in D^+$}.  Define $g=min\{|A|+1, |B|\}$, and $h= max\{|A|+1, |B|\}-g$.
 \begin{itemize}
 \item $|A|+1 > |B|$: In this case, on the $\lambda$ side, we remove $g=|B|$ beads from each runner, a total of $g^2$ moves.  We then allow the remaining  $h$ beads on runner $1$ to drop down $g$ places each, giving another $gh$ moves, for a total of $g(g+h)$.  The resulting core $\rho$ has $h>0$ beads on runner $1$ and is thus in $D^+$.

     The computation on the $\chi$ side is more difficult.  We first drop the $g=|B|$ beads on runner $1$ down one place each, filling up the space left because the part $1$ was removed by the involution $K_0$.  We remove $g$ pairs of beads, as total of $g^2$ moves as before, finally, we move the remaining $h-1$ (possibly $0$) beads on runner $p-1$ down, adding $g(h-1)$ moves. The total number of moves is $g(g+h)$ as before, the core $\rho'$ has $h-1$ beads on runner $p-1$, thus lying in $D^-$, and we have shown that in this case $\rho'=K_0(\rho)$.
 \item  $|A|+1 \le |B|$: In this case, we remove $g=|A|+1$ beads from each runner, a total of $g^2$ moves.  We then allow the remaining $h$  on runner $p-1$ to drop down $g$ places each, giving another $gh$ moves, for a total of $g(g+h)$.  The resulting core has $h$ (possibly 0) beads on runner $p-1$, but no beads on runner $1$, and is thus in $D^-$.

     On the $\chi$ side, we only have $g-1=|A|$ beads on runner $p-1$, so this is the maximum number of pairs which can be removed.  As before, we must first move $g-1$ beads down on runner $1$.  Thus removing the pairs requires a total of $g(g-1)$ moves.  There are $h+1$ beads left on runner $1$, suspended at a height of $g$, so an additional $g(h+1)$, bringing the total to $g(g+h)$ as required.  There are $h+1>0$ beads on runner $1$, so $\rho'$ lies in $D^+$, and we have again shown that $\rho'=K_0(\rho)$.

 \end{itemize}

  \textbf{Subcase 1.b:} \textsl{$\lambda \in D^-$}.  This is basically they same as Subcase 1.a, with the roles of $\lambda$ and $\chi$ reversed.

  Case 2:  $\lambda \in D$ has multiple parts on runner $0$.  Since none of the involutions affects runner $0$, we can remove the multiple parts from $\lambda$ and $\chi$, then apply Case 1.

\end{proof}

\section{Equivalences of extremal spin blocks.}

The Scopes involution $K_{i}$ $(0 \leq i \leq t)$ reduces blocks in the crystal graph. Kessar proved in certain cases  that external blocks of maximal $i$-strings in the block-reduced crystal graph are source algebra equivalent. Now we are about to consider what are the restrictions on $w$, in order that  $\rho^{w}$ be an external block in its $i$-string, for $0 \leq i \leq t$.

\begin{defn}

When $\rho^{w}$ is an external block in its maximal $i$-string, then we will say that the involution $K_{i}$   is a $w$-\textsl{allowed action}.

\end{defn}
\begin{flushleft}

\begin{lem}
Let $\rho$ be a $p$-core, with core $t$-tuple
\[c(\rho)=((\ell_{1},\epsilon_{1}),...,(\ell_{t},\epsilon_{t})),\].
\begin{enumerate}
	\item The involution $K_{0}$ is an $w$-allowed action if $w \leq \ell_{1}+\epsilon_1-1$
	\item The involution $K_{t}$ is an $w$-allowed action if $w \leq 2\ell_{t}+1$
	\item The involution $K_{i}$ ($1 < i <t$) is an $w$-allowed action if
		\begin{itemize}
			\item $w \leq (\ell_{i+1}-\ell_{i})\cdot(-1)^{\epsilon_{i}}$, for $\epsilon_{i}=\epsilon_{i+1}$
			\item $w \leq \ell_{i}+\ell_{i+1}$, for $\epsilon_{i}\neq \epsilon_{i+1}$
		\end{itemize}
\end{enumerate}
\end{lem}
\end{flushleft}

\begin{proof}

 The core $\rho$ is fixed. Let $v_{0}(i)$ be the weight $w$ of the  block on the last place in the block-reduced crystal graph that $\rho^{w}$ is external in direction $i$, and let $v_{1}(i)$ be $v_{0}(i)+1$, which is the weight $w$ of the block on the first place in the crystal graph that  $\rho^{w}$ is internal in direction $i$. By definition, only for external blocks do we get a $w$-allowed action, i.e. for blocks with weight $w$ satisfying $w \leq v_{0}(i)$. For every $w\geq v_1(i)$, $\rho^{w}$ will be internal because the crystal graph contains all of the translations of the exponents by positive integers for $\rho^{v_1(i)}$ and its adjacent blocks, \cite{AS}.

We assume $i$ fixed and write $v_1$ for $v_1(i)$. We investigate  $\rho^{v_{1}}$. Since  $\rho^{v_{1}}$ is internal with respect to $K_{i}$, there is a block $\mu^{v}$ such that  $\rho^{v_{1}}$ is obtained from $\mu^{v}$ by either restriction or induction, and is of rank either one less or, respectively, one greater. We now consider different cases of $i$, and different directions, rank-decreasing or rank-increasing.  Since, as proven in \cite{AS}, every block in the block-reduced crystal graph is the translation of a block of defect $0$, by the minimality of $v_1$, we see that the block which bounds it must be  of defect zero, so that $v$ is $0$.  Thus $\mu$ must be a core such that moving one bead to an adjacent runner produces the abacus representation of an element $\lambda$ of $RP_p$ which reduces in $v_1$ moves to the core $\rho$.  Since $\mu$, being a core, has no beads at all on the $p$ runner, the resulting element $\lambda$ of $RP_p$ is in fact a strict partition.

\begin{enumerate}
	\item $i=0$: Every removal of an $0$-good node from a core $\mu$ which will give $\rho^{v_1}$ requires that $\mu$ have the form $c(\mu)=((\ell_{1}',0),(\ell_{2},\epsilon_{2}),\dots,(\ell_{t},\epsilon_{t}))$.  The change which produces a restricted partition $\lambda$ in $\rho^{v_1}$ corresponds in terms of the abacus to removing the upper bead on the runner $1$ and putting it parallel in the $p$ runner. When we do so, this bead goes down $\ell_{1}'-1$ times and disappears, and $v_1 = \ell_{1}'-1$ (see \textsl{removing a $p$-bar}). Recalling that \[c(\rho)=((\ell_{1},\epsilon_{1}),...,(\ell_{t},\epsilon_{t})),\] we thus have  $\ell_{1}=\ell_1'-1$, so $v_1=\ell_{1}$ and $v_0 = \ell_1-1=\ell_{1}+{\epsilon_1}-1$. The block $\rho^{w}$ is internal if and only if $w\geq \ell_{1}$.

 The addition of an $0$-cogood node to the core $\mu$ which will give a partition in $\rho^{v_1}$ requires that $\mu$ have the form $c(\mu)=((\ell_{1}',1),(\ell_{2},\epsilon_{2}),\dots,(\ell_{t},\epsilon_{t}))$, with $\ell_1'-1 > 0$ (since otherwise the resulting partition is also a core, with weight $0$, and is not internal.
 The node which is added is the part $1$. This cancels the lowest bead on runner
 $p-1$, with gives one move, and the remaining $\ell_1'-1$ beads move down, so now we get $v_1=\ell_1'= \ell_1+1$, and $v_0= \ell_1$. Since in this case $\epsilon_1 = 1$, this gives the formula in the statement of the lemma.

	\item $i=t$: Every removal of an $t$-good node from a core $\mu^{w}$ for $c(\mu)=((\ell_{1},\epsilon_{1}),...,(\ell_{t}',1))$ corresponds in terms of the abacus to removing the upper bead on runner $t+1$ and putting it parallel in runner $t$. When we do so, this bead goes down $\ell_{t}'-1$ times, and then the two bottom-most beads on runners $t$ and $t+1$ are removed, all the beads on runner $t+1$ go down again, and the total number of moves is  by $v_1=2\ell_{t}'-2$. Let $\rho^{v_1}$ be the block obtained by removing one $t$-good node from the core $\mu$. It will satisfy $\ell_{t}=\ell_{t}'-2$, so $v_1=2\ell_{t}+2$ and thus $v_{0}=2\ell_{t}+1$.  The procedure for adding a $t$-cogood node is similar, except that it is the bottom bead which moves from the $t$ runner to the $t+1$ runner.  The total number of moves is the same.
	\item $0<i<t$: Every removal of an $i$-good node from a block $\mu^{v}$ for $c(\mu)=((\ell_{1},\epsilon_{1}),...,(\ell_{t},\epsilon_{t}))$ corresponds in terms of the abacus to removing a bead from runner $i+1$ or to removing a bead from runner $p-i$.
		\begin{itemize}
		\item First we examine the case that $\epsilon_{i}=\epsilon_{i+1}$:
		
				 \textsl{Case 1.} $\epsilon_{i}=\epsilon_{i+1}=0$: In this case, removing an $i$-good node corresponds to removing the upper bead in runner $i+1$ and putting it parallel in runner $i$. When we do so, this bead goes down till there is no empty places, i.e. $\ell_{i+1}'-\ell_{i}'-1$ moves. Let $\rho^{v_1}$ be the block after removing  an $i$-good node from the core $\mu$. The core t-tuple of $\rho$ satisfies $\ell_{i+1}=\ell_{i+1}'-1$ and $\ell_{i}=\ell_{i}'+1$, so $v_1=(\ell_{i+1}-\ell_{i}+1)$, and thus $v_{0}=\ell_{i+1}-\ell_{i}$.  The case of adding an $i$-cogood is similar.
 				
				 \textsl{Case 2.} $\epsilon_{i}=\epsilon_{i+1}=1$: In this case removing an $i$-good node is corresponds to removing the upper bead in runner $p-i$ and putting it parallel in runner $p-i-1$. When we do so, this bead goes down till there is no empty place i.e. $v_1=\ell_{i}'-\ell_{i+1}'-1$ times. Let $\rho^{v_1}$ be the block after removing  an $i$-good node from the core $\mu$. It satisfies $\ell_{i}=\ell_{i}'-1$ and $\ell_{i+1}=\ell_{i+1}'+1$, so $v_1=\ell_{i}-\ell_{i+1}+1$, and thus $v_{0}=\ell_{i}-\ell_{i+1}$.he case of adding an $i$-cogood is similar.
		\item Finally we examine the case that $\epsilon_{i}\neq \epsilon_{i+1}$:
				
		      Let $a_{1}p+(i+1)$ be the upper bead of runner $i+1$, and $a_{2}p+(p-i)$ be the upper bead of runner $p-i$.
		
		      \textsl{Case 1.} $a_{1}p+(i+1)>a_{2}p+(p-i)$: In this case removing an $i$-good node is corresponds to removing the upper bead from runner $i+1$, and put it parallel in runner $i$. When we do so, this bead goes down $\ell_{i+1}'-1$ times, then the bottom most bead of runner $p-i$ and the bead in runner $i$ are removed, and all the beads on runner $p-i$ go down again ($\ell_{i}'-1$ times).
		
		      \textsl{Case 2.} $a_{1}p+(i+1)\leq a_{2}p+(p-i)$: In this case removing an $i$-good node corresponds to removing the upper bead from runner $p-i$, and put it parallel in runner $p-i-1$. When we do so, this bead goes down $\ell_{i}'-1$ times, then the bottommost bead of runner $i+1$ and the bead in runner $p-i-1$ are removed, and all the beads on runner $i+1$ go down again ($\ell_{i+1}'-1$ times).
		
		      In both cases the weight is $v_1=\ell_{i}'+\ell_{i+1}'-1$. Let $\rho^{w}$ be the block after removing one time an $i$-good node from the core $\mu$. Its core t-tuple  satisfies $\ell_{i}=\ell_{i}'-1$ and $\ell_{i+1}=\ell_{i+1}'-1$, so $v_1=\ell_{i}+\ell_{i+1}+1$ and thus $v_{0}=\ell_{i}+\ell_{i+1}$.
		\end{itemize}
The procedure for adding a $i$-cogood node is is similar, except that the bead to be transferred is taken from the first place on the runner which is empty on the runner with which we are making the exchange.
\end{enumerate}
\end{proof}

\begin{defn} We say that two blocks $\rho^{w}$, $\sigma^{w}$ of the same weight are \textsl{allowed equivalent} if one can be obtained from the other by a sequence of $w$-allowed actions.
\end{defn}

\section{source algebra equivalences}

In this section we will demonstrate that a $w$-allowed action corresponds to a source algebra equivalence between appropriately chosen blocks labeled by the given cores.
We now define  $\mathcal M^{n-m}(\lambda,\chi)$ be the set of possible paths leading from the Young diagram of $\lambda$ to the Young diagram of $\chi$ by removing one node at a time in such a way that at every stage we have a strict partition.  If the Young diagram of $\chi$ is not a subset of $\lambda$, Then obviously the number will be zero.

For any core $\rho$, let $B_{\rho^w}$ be the block algebra of $\mathcal O\widetilde S_n$ with core $\rho$ and weight $w$ and let $B'_{\rho^w}$ be the corresponding block of $\mathcal O\widetilde A_n$ with core $\rho$ and weight $w$.

The key to proving the result we want will be the following definition.
Let $J_{n}$ be the set of strict partitions with core $\nu$ and weight $w$, and let $J_{m}$ be the set of strict partitions with core $\mu$ and weight $w$.

\begin{defn}\label{comp}\cite{KS}
 A\textsl{$w$-compatible pair}  $(\nu$, $\mu)$ for $i$ is defined to be a pair of cores such that:
\begin{enumerate}
\item  $K_{i}:J_{n} \rightarrow J_{m}$ is one-to-one and onto.
\item For any $\lambda\in J_{n}$ and $\chi \in J_{m}$,

\[
   |\mathcal M^{n-m}(\lambda,\chi)|= \left\{ \begin{array}{ll}
         0 & \mbox{if $\chi\neq K_{i}(\lambda)$};\\
        |\mathcal M^{n-m}(\nu,\mu)| & \mbox{if $\chi=K_{i}(\lambda)$}.\end{array} \right.
\]
\item $\epsilon(\lambda)+\epsilon(K_{i}(\lambda))=\epsilon(\mu)+\epsilon(\nu)$

\end{enumerate}
\end{defn}

We will now demonstrate that if $K_i$ gives a $w$-allowed action, the pair of cores $(\nu, \mu)$, is a $w$-compatible pair. At the end of the section, we will get the desired source algebra equivalence. We first prove some lemmas which will be needed to establish this result.

\begin{lem} Let $\nu$ be a core and let $i \in I$ be such that $\mu = K_i(\nu)$ is different from $\nu$ and of lower rank.  Then
\[
   |\mathcal M^{n-m}(\nu,\mu)|=\left\{
  \begin{array}{ll}
         (n-m)! & \mbox{if $i \neq 0$};\\
          \frac{(n-m)!}{2^{\frac{n-m-1}{2}}} & \mbox{if $i=0$}.
 \end{array} \right.
\]
\end{lem}

\begin{proof}

\begin{enumerate}
	\item $i\neq 0$: In this case the removal each of the $n-m$ $i$-boxes is independent of the removal of all the others, so we get  $(n-m)!$ possible orders in which we can remove these boxes.
	\item $i=0$: In this case in the Young diagram of $\mu$ the removable $0$-boxes comes in pairs except the last row that has only a single 0. Thus $n-m$ is odd, and we have $\frac{n-m-1}{2}$ pair of 0-boxes. The number of possible orders in which we can  remove these boxes without restriction is $(n-m)!$. However, each adjacent pair we must remove  first the extremal node and then the one next to it, a consideration which divides the number of permissible orderings by 2.   For each adjacent pair we must divide by $2$, so altogether we have to divide by $2^{\frac{n-m-1}{2}}$.  This reduces the number of permissible ordering to $\frac{(n-m)!}{2^{\frac{n-m-1}{2}}}$.
\end{enumerate}

\end{proof}

As in Definition \ref{comp}, let $J_{n}$ be the set of strict partitions with core $\nu$ and weight $w$, and let $J_{m}$ be the set of strict partitions with core $\mu$ and weight $w$. Take $\lambda\in J_{n}$, $\chi\in J_{m}$, where $n>m$ by our assumptions on $\nu$ and $\mu$.  To go back from $\nu$ to $\lambda$ (or from $\mu$ to $\chi$ ) we have to do $w$ moves which correspond to adding $p$-bars. There are three kinds of moves:
\begin{itemize}
	\item Moving a bead up on its runner.
	\item Inserting a pair of beads to runners $i$ and $p-i$ where the bottom place in each runner is empty.
    \item creating a bead on runner $0$.

\end{itemize}

\begin{lem} Let $\nu$ be a core. The number of actions that needed to insert $n$ pairs of beads on runners $i,p-i$ for some $i>0$ in the abacus representation of $\nu$ is:

	 $$n^{2}+\ell_{i}n$$

\end{lem}
\begin{proof} By induction.

	For $n=1$: To insert one pair we have to lift the $\ell_{i}$ beads on runner $i$ or on runner $p-i$ upwards, that is, $\ell_{i}$ moves and then insert the pair, that is, one more move. Altogether, we make $\ell_{i}+1$ moves.
	
	Assume that the lemma is true for $n$, i.e., the number of actions that needed to insert $n$ pairs of beads on runners $i,p-i$ in the abacus representation of $\nu$ is $n^{2}+\ell_{i}n$. Now we prove that lemma is true for $n+1$:

To insert a pair after inserting $n$ pairs: First we have to lift the $n+\ell_{i}$ beads that are on runner $i$ or on runner $p-i$ upwards, that is $n+\ell_{i}$ actions, Second, we have to lift the $n$ beads that are on the other runner upwards, that is, another $n$ moves. We can then insert the new pair, another move, altogether doing $2n+\ell_{i}+1$ new moves. Inserting $n$ pairs, by the  assumption   is $n^{2}+n\ell_{i}$ moves, so we have $$n^{2}+n\ell_{i}+2n+\ell_{i}+1=(n+1)^2+\ell_{i}(n+1).$$

\end{proof}

\begin{corl}Suppose that $K_i$ is a $w$-allowed action.

 It is not possible to insert a pair of beads on runners $j$ and $p-j$ in the following cases:
\begin{enumerate}
	\item $j=i$ for $0<i<t$, $\epsilon_{i}=\epsilon_{i+1}=1$, and  $\ell_i>\ell_{i+1}$.
	\item $j=i$ for $0<i<t$ , $\epsilon_{i}=\epsilon_{i+1}=0$ and $2\ell_{i}+1\geq \ell_{i+1}>\ell_{i}$.
	\item $i=0$ and $j=1$.
\end{enumerate}
	Nor is it possible to insert more then one pair of beads on runners $t$ and $t+1$ when $\ell_{t}>0$.

\end{corl}
\begin{proof} We divide into cases:

\begin{itemize}
	\item $i \neq t$:
\begin{enumerate}
	\item In this case, we know that $w\leq \ell_{i}-\ell_{i+1}$ and inserting this pair is $\ell_{i}+1$ moves, i.e. $\ell_{i}+1 \leq w \leq
			\ell_{i}-\ell_{i+1}$, which implies $\ell_{i+1} \leq -1$, a contradiction to the definition of the $\ell_i$.
	\item In this case we know that $w\leq \ell_{i+1}-\ell_{i}$ and inserting this pair is $\ell_{i}+1$ actions, i.e. $\ell_{i}+1 \leq w \leq
			\ell_{i+1}-\ell_{i}\Rightarrow 2\ell_{i}+1 \leq \ell_{i+1}$.
	\item In this case we know that $w\leq \ell_{1}-1$ if $\epsilon_{1}=0$ and $w\leq \ell_{1}$ if $\epsilon_{1}=1$. Inserting this pair is $\ell_{1}+1$ actions, i.e. $\ell_{1}+1 \leq w \leq 		
			\ell_{1} \Rightarrow 1\leq 0 $, a contradiction.
\end{enumerate}
	\item $i = t$: In this case we know that $w\leq 2\ell_{t}+1$ and inserting more then one pair is at least $2\ell_{t}+4$ actions, a contradiction.
\end{itemize}

\end{proof}

\begin{defn} We say that a bead \textsl{can be reduced} if and only if the bead can move from runner $i$ to an empty place at the same height, in runner $i-1$ for $2<i<p-i$, or from runner 1 to an empty place at one height less, in runner $p-1$ for height bigger than 1, or from runner 1 at height 1 to disappearance of the bead. The move from runner $1$ to runner $p-1$ is actually made in two steps, first to runner $0$ and then to runner $p-1$.  If the first move is blocked by a bead on runner $0$, then we first move the bead on runner $0$ to runner $p-1$, and only afterwards move the bead on runner $1$ to runner $0$. This action will be called \textsl{to reduce a bead}. Such a reduction reduces the rank by $1$ or $2$, the latter only when moving the bead from runner $1$ to runner $p-1$.
\end{defn}

\begin{remark}
The number $|\mathcal M^{n-m}(\nu,\mu)|$, in terms of the abacus, is the number of possible orders in which to reduce beads and the cases in which this number is positive are:
\begin{enumerate}
	\item $K_i$, $0<i<t$:
\begin{itemize}

	\item $\epsilon_{i}=1$, $\epsilon_{i+1}=0$,$\ell_i+\ell_{i+1}>0$.
	\item $\epsilon_{i}=\epsilon_{i+1}=0$, $\ell_{i+1}>\ell_{i}$.
	\item $\epsilon_{i}=\epsilon_{i+1}=1$, $\ell_{i+1}<\ell_{i}$.
\end{itemize}
	\item $K_0$ ,$\ell_{1}>0$  $\epsilon_{1}=0$.
	\item $K_t$ , $\ell_{t}>0$, $\epsilon_{t}=1$.
\end{enumerate}
The complementary cases all increase the rank.  There are no $w$-allowed actions which leave the rank fixed, except those which are trivial because there are identical configurations of beads on each pair of interchanged runners. The action $K_0$ is never trivial.

\end{remark}

We now have enough information to prove the main result of this section:

\begin{prop} If $K_i$ gives a $w$-allowed action with respect to a pair of cores
$\nu$, $\mu$, then ($\nu$, $\mu$) is a $w$-compatible pair via $K_i$.
\end{prop}

\begin{proof}
\begin{enumerate}
	\item $K_i$ is an involution. By Lemma 2.1 we see that for every strict partition $\lambda \in J_n$ there is a suitable strict partition $\chi \in J_m$ obtained by $K_i$.
	\item  We need to compare $|\mathcal M^{n-m}(\lambda,\chi)|$ and $|\mathcal M^{n-m}(\nu,\mu)|$.  The strict partition $\lambda$ is obtained from the core $\nu$ by adding $w$ $p$-bars.  To the extent that these moves take place on runners not affected by $K_i$, the same moves will be involved in obtaining $\chi$ from $\mu$.  Thus the only moves which affect the possibility of reducing beads are on the runners affected by $K_i$, or on the runner $0$ in the case of $i=1$.  We consider the cases listed in the remark, which are the only cases relevant when $K_i$ reduces the rank.
\begin{enumerate}

	\item $K_i$, $0<i<t$:
\begin{itemize}

	\item \underline{Case 1: $\epsilon_{i}=1$, $\epsilon_{i+1}=0$, $\ell_i+\ell_{i+1}>0$.} Moving beads up does not block any bead reductions. By Corollary 4.2.1, we cannot add pairs unless $\ell_i >\ell_{i+1}$. If we add $d>0$ pairs to $\nu$, then we must use up $d\ell_i+d^2$ moves, and this is greater than or equal to $\ell_i+d$ so
$$\ell_i+d \leq w \leq \ell_i+\ell_{i+1}.$$
The number of moves remaining is then less than or equal to
$$ \ell_{i+1}-d.$$
Thus in the remaining moves on runner $i$, the topmost bead cannot rise above the topmost bead in runner $i+1$.  The number of bead on runner $i+1$ which can be reduced decreases by $d$, but the number of beads which can be reduced on runner $i$ increases by the same $d$, so the total number of beads which can be reduced is the same as in $\nu$.

	\item \underline{Case 2: $\epsilon_{i}=\epsilon_{i+1}=0$, $\ell_{i+1}>\ell_{i}$.} Since $w$ is the positive difference, we cannot move any beads on runner $i$ up past the topmost bead of runner $i+1$.  We can add pairs only to runner $i$ and  only when  $$d \leq d\ell_i+d^2 \leq w \leq\ell_{i+1}.$$  The beads whose reduction is blocked on runner $i+1$ will be compensated for by $d$ new beads on runner $i$, so again the total number will be the same.
	\item \underline{Case 3:$\epsilon_{i}=\epsilon_{i+1}=1$}, $\ell_{i+1}<\ell_{i}$.  This case is similar to Case 2, with $i$ and $i+1$ reversed.
\end{itemize}
	\item \underline{$K_0$ ,$\ell_{1}>0$  $\epsilon_{1}=0$.}
In this case, there are no pairs which can be added, and no beads on runner $p-1$ to block the reduction of beads from runner $1$. However, here we must deal with the problem of beads on runner $0$, in the manner described above.
	\item \underline{$K_t$ , $\ell_{t}>0$, $\epsilon_{t}=1$.} In this case, just moving beads up will not block any reductions.  If a pair is added, it is only one, and adding it uses up $\ell_t +1$ moves, while $w \leq 2\ell_t+1$, so the new bead on runner $t$ cannot rise above the topmost bead of runner $t+1$.
\end{enumerate}

	\item  As we said before, $\lambda$ belongs to the blocks $\nu^{w}$, and $\chi$ to the block $\mu^{w}$. For every core $\rho$, we defined   $\epsilon(\rho)\equiv\left|\rho\right|+n(\rho)$(mod 2). Adding a $p$-bar to the core $\rho$, can be done in one of the two ways: One way is lifting up one bead on its runner, so the number of parts does not change but $\left|\rho\right|$ becomes $\left|\rho\right|+p$. The other way is adding a pair of beads ($i$, $p-i$), so $n(\rho)$ becomes $n(\rho)+2$ and $\left|\rho\right|$ becomes $\left|\rho\right|+p$, (recall $p$ is odd). In both ways the parity of the partition is changed from even to odd, or from odd to even. For $w$ additions of $p$ we can summarize the result: If $w$ is even, the parity changes $w$ times, and therefore $\epsilon(\rho^{w})=\epsilon(\rho)$, and if $w$ is odd, then the parity changes according to the following scheme: $\epsilon(\rho^{w})=1-\epsilon(\rho)$. Therefore, for an even $w$ we get $\epsilon(\lambda)+\epsilon(\chi)=\epsilon(\nu)+\epsilon(\mu)$ and for odd $w$ we get $\epsilon(\lambda)=1-\epsilon(\nu)$ and $\epsilon(\chi)=1-\epsilon(\mu)$ so, $\epsilon(\lambda)+\epsilon(\chi)=1-\epsilon(\nu)+1-\epsilon(\mu)\equiv\epsilon(\nu)+\epsilon(\mu) (mod 2)$.

\end{enumerate}
\end{proof}

We now come to the source algebra equivalence and the problem of the crossovers. We will work over a modular system $(\mathcal K, \mathcal O, \mathcal F)$.

 Let $\nu$ be a core of rank $n-wp$ and let $\mu$ be a core of rank $m-wp$, with $n>m$. Let $b \in \widetilde A_n$ represent the block idempotent of $B_{\nu^w}$ or $B'_{\nu^w}$, and let $c \in \widetilde A_m$ represent the block idempotent of $B_{\mu^w}$ or $B'_{\mu^w}$.
If $K_i$ is a $w$-allowed action with $\mu = K_i(\nu)$,
then we are trying to prove that one of the following holds:
\begin{enumerate}
\item If the parities are the same, $B_{\nu^w}$ is source algebra equivalent to $B_{\mu^w}$, and $B'_{\nu^w}$ is source algebra equivalent to $B'_{\mu^w}$,
\item If the parities are different, $B_{\mu^w}$ is source algebra equivalent to $B'_{\mu^w}$, and $B'_{\nu^w}$ is source algebra equivalent to $B_{\mu^w}$,
\end{enumerate}

 The first set of equivalences was essentially proven in \cite{K}.  However, since the main thrust of that paper was the Donovan conjecture, it is rather hard to extract the particular result that we need.  Therefore, we recast our results in a form which will allow us to apply Theorem 2.5 in \cite{HK}. For the second set of equivalences, we will cite \cite{KS}, which used permutation modules.

 For any strict partition $\lambda$, let $\theta_\lambda^\pm$ or $\eta_\lambda^\pm$ be the corresponding irreducible character or characters of $\widetilde S_n$  or $\widetilde A_n$, depending on the parity.  In $\widetilde S_n$, if $\epsilon(\lambda)= 0$ we get one character $\theta_\lambda^+$, and if $\epsilon(\lambda)= 1$ we get two associate characters $\theta_\lambda^\pm$.   In $\widetilde A_n$, if $\epsilon(\lambda)= 1$ we get one character $\eta_\lambda^+$, and if $\epsilon(\lambda)= 0$ we get two conjugate characters $\eta_\lambda^\pm$.
 \begin{defn}We define
 $r(\theta_\lambda^\pm,\theta_\chi^\pm,bc)$ to be the
  number of constituents of one of the characters $\theta_\lambda^\pm$ in the block with idempotent $bc$ after inducing of one of the characters $\theta_\chi^\pm$ of $\widetilde S_m$ to $\widetilde S_n$. The corresponding number for the $\widetilde A_n$ will be denoted by $r'$.   \end{defn}

 We are restricting ourselves henceforward to the the case where the parities are the same, so either both characters belong to associate pairs or both are self-associate.  This means that the same formula holds as well for restriction from $\widetilde A_n$ to $\widetilde A_m$, by standard Clifford theory.

 In order to quote Theorem 2.5 from \cite{HK}, we must make a few general definitions.  Let $G$ be a group and $H$ a subgroup containing a $p$-group $D$.   Set $A= \mathcal OG$ and let
 $$Br_D: A^D \rightarrow \mathcal F C_G(D)$$ be the Brauer homomorphism.  For any idempotent $u$ of $(\mathcal OG)^H$, we let $m_{H,D}(u\mathcal OGu)$ be the number of idempotents $i$ in a primitive idempotent decomposition of $u$ in $(u\mathcal OGu)^H$ for which $Br_D(i)$ is non-zero.

 As in the case to which we wish to apply the theorem, we let $b$ be the idempotent of a block of $\mathcal O G$ and let $c$ be the idempotent of a block of $\mathcal OH$, assuming that the two blocks have a common defect group $D$.

\begin{lem}\label{K}  Let  $\alpha := n-m $ and let
$\beta := |{\mathcal M}^\a (\nu, \mu)|$.   Suppose that   $ \nu $ and $ \mu $ form
a  $w$-compatible pair and  that $\nu $ and $\mu $ have the same parities.
Let $\epsilon(\a)$ be $0$ or $1$ as $\a$ is even or odd.
Then   $|\Irr (\tilde S_n, b)| =   |\Irr (\tilde S_m, c)| $
and  for each
$\theta_\lambda^\pm \in  \Irr(\tilde S_n, b )$,
$$ \sum_{\theta_\tau^\pm \in \Irr(\widetilde S_m, c)}r(\theta_\lambda^\pm, \theta_\tau^\pm, bc )=
2^{\frac{\alpha -\epsilon(\a)}{2}}\beta $$ and for each
$\theta_\chi^\pm \in  \Irr (\tilde S_m, c )$,
$$ \sum_{\theta_\tau^\pm \in \Irr (\widetilde S_n, b)}r(\theta_\tau^\pm, \theta_\chi^\pm, bc )=
2^{\frac{\alpha -\epsilon(\a)}{2}}\beta . $$

\end{lem}

\begin{proof}  Let $\lambda $   be a  strict partition of $n$   and
$\chi $  a strict partition of $m$. Let $\theta_\lambda^\pm $
be an irreducible  character of $\tilde S_n $  corresponding  to
$\lambda $
and $\theta_\chi^\pm $ an irreducible  character of
$\tilde  S_m $ corresponding  to
$\chi $.
It follows from the branching rules
(see  for example \cite{St}) that if $\theta_\lambda^\pm $ is a  constituent of
$Ind_{\tilde  S_m}^{\tilde S_n} (\theta_\chi^\pm)$, then
${\mathcal M}^\a(\lambda, \chi) $ is non-empty, which implies that
${\mathcal M}^\a(\nu, \mu) $ is also non-empty.
\begin{enumerate}
\item $\a \neq 1$.
Since  in our situation,
$\nu $ and $\mu $ have the same parities, then $\a$ is odd if $i=0$ and $\a$ is even if $i>0$, for $i$ such that $K_i(\mu)=\nu$.
We can calculate the coefficients of induced characters using the branching rules in \cite{St}.

 If $\a$ is even, then $i \neq 0$, so that
$\lambda $ and $\chi $ have the same  number of parts, and
the multiplicity of  $\theta_\lambda^\pm $ as a  constituent of
$Ind_{\tilde  S_m}^{\tilde S_n}(\theta_\chi^\pm)$    is
$2^{\frac{\alpha}{2}}\beta$
if  $\epsilon(\chi)=0 $  and is
$2^{\frac{\alpha}{2}-1}\beta$ if
$\epsilon(\chi)=1 $.

If $\alpha $ is odd, then $i=0$, so that $\lambda$ has one more part then $\chi$.  If $\a > 1$,
the multiplicity of  $\theta_\lambda^\pm $ as a  constituent of
$Ind_{\tilde  S_m}^{\tilde S_n} (\theta_\chi^\pm)$ is $2^{\frac{\alpha-1}{2}}\beta$
if $\epsilon(\chi) = 0$ and is $2^{\frac{\alpha-3}{2}}\beta$ if $\epsilon(\chi) = 1$.

 We can summarize these branching rules in a single formula, recalling that $\epsilon(\a)$ is the parity of $\a$.  Let $\beta =  |\mathcal M^\a(\nu,\mu)|$.
 Then, if  $\a \neq 1$,

 $$r(\theta_\lambda^\pm,\theta_\chi^\pm,bc)
 =2^{\frac{\a-\epsilon(\a)-\epsilon(\lambda)-\epsilon(\chi)}{2}}\mathcal \beta.$$

 We now verify that the sums in the statement of the lemma are correct. When $\epsilon(\chi) = 0$, there is only one character in each of the sums in the statement of the lemma with non-zero coefficients, and that coefficient has the correct value
 required for the lemma.  When $\epsilon(\chi) = 1$, then there are exactly two with non-zero coefficient, and the sum of those two coefficients has the correct value.

\item $\a=1$
  We have $\beta=1$, since the only way to have the parities preserved is by having $\nu$ obtained from $\mu$ by adding the part $1$. In this case  $\theta_\chi^+$  lifts to a unique character, either
$\theta_\lambda^+ $ or $\theta_\lambda^- $, with multiplicity $1$ (and similarly for  $\theta_\chi^-$, when $\epsilon(\chi) = 1$.)

 In the special case that $\a =1$ and $\epsilon(\chi) = 0$,
  $$r(\theta_\lambda^+,\theta_\chi^+,bc)=1,$$ giving both of the sums in the lemma (recalling that $\beta = 1$ and $\a - \epsilon(\a) = 0$.

 In the special case that $\a =1$ and $\epsilon(\chi) = 1$, then there is a specific correspondence between the paired characters, so
$$r(\theta_\lambda^+,\theta_\chi^+,bc)+r(\theta_\lambda^+,\theta_\chi^-,bc)=1,$$ and
 $$r(\theta_\lambda^-,\theta_\chi^+,bc)+r(\theta_\lambda^-,\theta_\chi^-,bc)=1,$$
which give the first sum in the statement of the lemma.

 It is equally true that
  $$r(\theta_\lambda^+,\theta_\chi^+,bc)+r(\theta_\lambda^-,\theta_\chi^+,bc)=1,$$ and
 $$r(\theta_\lambda^+,\theta_\chi^-,bc)+r(\theta_\lambda^-,\theta_\chi^-,bc)=1,$$ which gives the second sum in the statement of the lemma.

\end{enumerate}

The similar results for $\tilde A_n $ and $\tilde A_m $  can be proven in an almost identical fashion.  It is done in full in \cite{L2}.
\end{proof}

 \begin{thm*} \cite[Theorem 2.5]{HK} Let $u$ be an idempotent of $(c\mathcal OGbc)$.  Then for any character $\phi \in Irr(G,b)$
 $$\sum_{\psi \in Irr(H,c)} r(\phi,\psi,u) \geq m_{H,D}(u\mathcal OGu),$$
 and for any character $\psi \in Irr(H,c)$ we have
 $$\sum_{\phi \in Irr(G,b)} r(\phi,\psi,u) \geq m_{H,D}(u\mathcal OGu).$$
     Further, if $|Irr(H,c)| = |Irr(G,b)|$ and if
 $$\sum_{\psi \in Irr(H,c),\phi \in Irr(G,b)} r(\phi,\psi,u) \leq m_{H,D}(u\mathcal OGu)|Irr(H,c)|,$$
 then for any primitive idempotent $i$ of $(u\mathcal OGu)$, the image $Br_D$(i) is non-zero and the map $\mathcal OHc \rightarrow i\mathcal OGi$ given by $x\mapsto ix$ is an isomorphism of interior $H$-algebras.  In particular, if the above holds, then $\mathcal OGb$ and $\mathcal OHc$ are Puig equivalent.
 \end{thm*}

Since we are assuming that $(\nu, \mu)$ is a $w$-compatible pair, we are in the situation of the second part of the theorem. We have already calculated the sum in Lemma \ref{K}.
 To complete the proof, we need the following modified version of an unpublished lemma by Kessar from an early version of \cite{KS}.

 \begin{lem}
 With the notation above,
 \begin{itemize}
\item If $\a$ is odd, $m_{H,D}(\mathcal OGbc)=2^{\frac{\alpha -1}{2}}\beta $.
\item If $\a$ is even,$m_{H,D}(\mathcal OGbc)=2^{\frac{\alpha}{2}}\beta $.
\end{itemize}
 \end{lem}
\begin{proof}

 By Brauer's Main Theorem (see, for example, \cite{A}), the image $Br_D(b)$ of the block idempotent $b$ in $\mathcal FC_G(D)$ is the $N_G(D)$ conjugacy class sum of block idempotents of blocks whose image under the projection $\pi:\mathcal FC_G(D)\rightarrow \mathcal FC_G(D)/Z(D)$ is of defect zero.  Let $\bar b$ be $\pi \circ Br_D$, and define $\bar c$ in a similar fashion for $H$.

In order to apply Theorem 2.5 of  \cite{HK} in the proof of our main theorem, we must first calculate
the quantity $ m_{H, D}(c {\mathcal O} G bc)$, which is the number of primitive idempotents in an idempotent decomposition of $(c \mathcal O Gbc)^H$ which are non-zero under the action of $\pi \circ Br_D$.   By Proposition 2.6 of \cite{HK},
this is the same as  the number of idempotents  in  a primitive idempotent
decomposition
of  the algebra $ (\bar c  \mathcal F(C_G(D)/Z(D)\bar b \bar c)^{N_H(D)}$, where
$N_H(D)$ acts on  $\mathcal FC_G(D)$  through the inclusion of $H$ in $G$.  This number is the product of the degree of the defect zero block, times the length of the orbit under the action of $N_H(D)$.

According to  the local structure of blocks of the double covers of
the  symmetric groups given in Cabanes  \cite {Ca}, the defect group is
isomorphic to the defect group of $\widetilde S_{|pw|}$.
There is  a subgroup, $\widetilde S_{|pw|}\widetilde S_{|\nu|} $   of
$\widetilde S_n $  lifting a  Young subgroup of $S_n $ isomorphic to
$S_{|pw|}\times  S_{|\nu|} $  such that
$C_G(D)/Z(D)  \cong  \widetilde S_{|\nu|} $.  The inclusion of $H$ in $G$ and  the
defect group $D$
can be chosen such
that  $C_H(D)/Z(D) \cong
\widetilde S_{|\mu|} $,  where the induced
embedding of
$\widetilde S_{|\mu| }$    in $\widetilde S_{|\nu|}$ is the standard embedding,
 and  $N_H(D) \cong   \widetilde  S_{|\mu|}
N_{\widetilde S_{|pw|}}(D) $.  Let  $\Lambda (b)$ be the  set of
ordinary irreducible characters  of $\widetilde S_{|\nu|} $ corresponding to the
partition $\nu $, which is either $\{\theta_\nu^+\}$ or $\{\theta_\nu^+,\theta_\nu^-\}$, depending on parity.  Similarly, let
 $\Lambda (c)$ is the  set of
ordinary irreducible characters  of $\widetilde S_{|\mu|} $ corresponding to the
partition $\mu $.

The groups $\widetilde S_{pw}$ and $\widetilde S_{|\nu|} $ do not centralize each other, even though the sets on which their images in $\widetilde S_{n}$ act are disjoint, because commuting transpositions multiplies the product by the central involution $z$ (which acts in spin representations as $-1$.) However, the group $\widetilde A_{|pw|}$  centralizes $\widetilde  S_{|\nu|} $ since all its elements are even products of transpositions, hence
$$ (\bar c  \mathcal F(C_G(D)/Z(D))\bar b \bar c)^{N_H(D)} =
(\bar c \mathcal F\widetilde S_{|\nu|} \bar b \bar c) ^{<\widetilde S_{|\mu|},\sigma>}, $$
where $\sigma$  in $\widetilde S_{|pw|}-\tilde A_{|pw|}$ is the lifting of a transposition, so that $\sigma^2$ is central.
 Also, since $\sigma $
normalizes  $\widetilde S_{|\mu|}$, $\sigma $ acts on
$(\bar c \mathcal F \widetilde S_{|\nu|} \bar b \bar c) ^{ \tilde S_{|\mu|}} $ and thus
we can compute $$ (\bar c \mathcal F (C_G(D)/Z(D))\bar b \bar c)^{N_H(D)} $$ by first
computing the algebra $(\bar c \mathcal F\widetilde S_{|\nu|} \bar b \bar c) ^{ \widetilde S_{|\mu|}}$ and
then taking  fixed points under $\sigma$.

 The bimodule
$W=(\bar c \mathcal F\widetilde S_{|\nu|} \bar b \bar c) ^{\widetilde S_{|\mu|}} $  is
isomorphic to   $$ End_{\mathcal F (\widetilde S_{|\nu|}\times  \widetilde S_{|\mu|}^{op})}
(\mathcal F\tilde S_{|\nu|}\bar b \bar c),$$ the  algebra of
$\mathcal F(\widetilde S_{|\nu|} \times \widetilde S_{|\mu|}^{op})$  invariant
endomorphisms
of the $\mathcal F(\tilde S_{|\nu|} \times \tilde S_{|\mu|}^{op})$-module
$\mathcal F\tilde S_{|\nu|}\bar b \bar c $.   The map taking $\bar b \bar c$ to $w \in W$ is obviously a homomorphism of left modules, and it is a homomorphism for the right action as well because $w$ is fixed under conjugation by elements of $\widetilde S_{|\mu|}$.

 Since $\bar b $ and $\bar c$ are
of defect zero,  the
$\mathcal F(\tilde S_{\nu} \times \tilde S_{\mu}^{op})$-module
$\mathcal F\tilde S_{|\nu|}\bar b \bar c $ is isomorphic to
$$ \sum_{\phi \in \Lambda(b), \psi \in \Lambda(c)} d_{\phi, \psi}
V_{\phi} \otimes V_{\psi}, $$ where $V_{\phi}$  and   $V_{\psi}$ are  simple
projective modules for  $ \mathcal F\widetilde S_{|\nu|}$ and $\mathcal F\widetilde S_{|\mu|} $
corresponding to the ordinary irreducible characters $\phi $ and $\psi $
respectively, and where
$d_{\phi, \psi}$ is  the multiplicity of $\phi $ in
$ Ind_{\widetilde S_{|\mu|}}^{\widetilde  S_{|\nu|}} (\psi) $.
Thus $(\bar c \mathcal F \widetilde S_{|\nu|} \bar b \bar c) ^{ \widetilde S_{|\mu|}} $ is isomorphic to the semi-simple algebra
$$\prod_{\phi \in \Lambda(b), \psi \in \Lambda(c)} Mat_{d_{\phi,\psi}}(\mathcal F).$$
\begin{enumerate}
\item $\a$ is odd.
Let us consider the case  $\epsilon(\nu)=0 $, $\epsilon (\mu)=0 $.
In this case  $\Lambda(b)$   consists of the unique irreducible character
$\theta_\nu^+ $ of
$\widetilde S_{|\nu|}$ corresponding to $\nu $ and
$\Lambda(c)$   consists of the unique irreducible character $\theta_\mu^+ $, of
$\widetilde S_{|\mu|}$ corresponding to $\mu $. The multiplicity  of $\theta_\nu^+$
in $Ind_{\widetilde S_{|\mu|}}^{\tilde  S_{|\nu|}} (\theta_\mu^+)  $ is
$2^{\frac{\alpha -1}{2}}\beta $, so
$(\bar c \mathcal F\widetilde S_{|\nu|} \bar b \bar c) ^{ \widetilde S_{|\mu|}} $ is a matrix
algebra of size $2^{\frac{\alpha -1}{2}}\beta $.

 In particular,  $\sigma$ acts as an inner automorphism
 on
$(\bar c k\tilde S_{|\nu|} \bar b \bar c) ^{ \tilde S_{|\mu|} }$.  Since
$\sigma ^2 $ is central and thus acts  as the identity,  and since  $p$ is of odd
characteristic, we may assume that the  action
of $\sigma$ is through a diagonal matrix with $1$'s and $-1$'s on the
diagonal. Thus the fixed points of this action  are block diagonal matrices corresponding to the decomposition into eigenspaces of $\sigma$.  The total number of primitive idempotents remains the same, equal to the total degree of the block diagonal matrix. It
follows that  the number of idempotents  in  any  primitive idempotent
decomposition of
$(\bar c \mathcal F\widetilde S_{|\nu|} \bar b \bar c) ^{ <\widetilde S_{|\mu|}, \sigma>}$ is
$2^{\frac{\alpha -1}{2}}\beta $.

Now let us consider the case  $\epsilon(\nu)=1 $, $\epsilon (\mu)=1 $

Then each of  $\Lambda (b)$ and  $\Lambda (c)$ consist of two characters.
\begin{enumerate}
\item $\a > 1$:
The multiplicity of any  irreducible character in $\Lambda (b)$ in the
induced character
of any  irreducible character  in
$\Lambda (c)$ is $2^{\frac{\alpha -3}{2}}\beta $.
 Thus,
$(\bar c \mathcal F\widetilde S_{|\nu|} \bar b \bar c) ^{ \widetilde S_{|\mu|}} $ is a direct
product of four matrix algebras each of size $2^{\frac{\alpha -3}{2}}\beta $.
The element $\sigma $ permutes  these  matrix factors in pairs, so
$(\bar c \mathcal F\widetilde S_{|\nu|} \bar b \bar c) ^{ <\widetilde S_{|\mu|}, \sigma>}$
is isomorphic to a direct  product of two matrix algebras  each of size
 $2^{\frac{\alpha -3}{2}}\beta $. Thus the number of idempotents  in  a
primitive idempotent
decomposition of
$(\bar c \mathcal F\widetilde S_{|\nu|} \bar b \bar c) ^{ <\widetilde S_{|\mu|}, \sigma>}$ is
$2^{\frac{\alpha -1}{2}}\beta $.
\item $\a =1$: In this case $\beta=1$, and there is a pairing between the elements of $\Lambda(b)$ and $\Lambda(c)$, so that the number of constituents is either $0$ or $1$.  We may assume that in this special case, $\theta_\mu^+$ lifts to $\theta_\nu^+$ and $\theta_\mu^-$ lifts to $\theta_\nu^-$.
The total number of idempotents lifting one of the elements of $\Lambda(c)$ is $1$, but this is exactly equal to  $2^{\frac{\alpha -1}{2}}\beta $, as in the case $\a > 1$.
\end{enumerate}
 Thus the number of idempotents  in  a
primitive idempotent
decomposition of
$(\bar c \mathcal F\widetilde S_{|\nu|} \bar b \bar c) ^{ <\widetilde S_{|\mu|}, \sigma >}$ is
$2^{\frac{\alpha -1}{2}}\beta $.

\item $\a$ is even.

Let us consider the case  $\epsilon(\nu)=0 $, $\epsilon (\mu)=0 $.
In this case  $\Lambda(b)$   consists of the unique irreducible character
$\theta_\nu^+ $ of
$\tilde S_{|\nu|}$ corresponding to $\nu $ and
$\Lambda(c)$   consists of the unique irreducible character $\theta_\mu^+ $, of
$\tilde S_{|\mu|}$ corresponding to $\mu $. The multiplicity  of $\theta_\nu^+  $
in $Ind_{\tilde S_{|\mu|}}^{\tilde  S_{|\nu|}} (\theta_\mu^+)  $ is
$2^{\frac{\alpha }{2}}\beta $, so
$(\bar c k\tilde S_{|\nu|} \bar b \bar c) ^{ \tilde S_{|\mu|}} $ is a matrix
algebra of size $2^{\frac{\alpha }{2}}\beta $.

Now lets consider the case  $\epsilon(\nu)=1 $, $\epsilon (\mu)=1 $.
Then each of  $\Lambda (b)$ and  $\Lambda (c)$ consist of two characters,
and the multiplicity of any  irreducible character in $\Lambda (b)$ in the
induced character
of any  irreducible character  in
$\Lambda (c)$ is $2^{\frac{\alpha }{2}-1}\beta $. Thus,
$(\bar c \mathcal F\widetilde S_{|\nu|} \bar b \bar c) ^{ \widetilde S_{|\mu|}} $ is a direct
product of four matrix algebras each of size $2^{\frac{\alpha }{2}-1}\beta $.
The element $\sigma $ permutes  these  matrix factors in pairs, so
$(\bar c \mathcal F\widetilde S_{|\nu|} \bar b \bar c) ^{ <\widetilde S_{|\mu|}, \sigma>}$
is isomorphic to a direct  product of two matrix algebras  each of size
 $2^{\frac{\alpha}{2}-1}\beta $. Thus the number of idempotents  in  a
primitive idempotent
decomposition of
$(\bar c \mathcal F\tilde S_{|\nu|} \bar b \bar c) ^{ <\tilde S_{|\mu|}, \sigma >}$ is
$2^{\frac{\alpha }{2}}\beta $.

\end{enumerate}
\end{proof}

We remind the reader of the notation $B_{\rho^w}$ and $B'_{\rho^w}$ given in the introduction for blocks of $\wt S_n$ and $\wt A_n$.

\begin{thm} Suppose $\nu^w$ and $\mu^w$ are extremal points of an $i$-string in the block-reduced crystal graph.
\begin{enumerate}
\item
If the parities are the same, $B_{\nu^w}$ is source algebra equivalent to $B_{\mu^w}$, and $B'_{\nu^w}$ is source algebra equivalent to $B'_{\mu^w}$.
\item If the parities are  different,$B_{\nu^w}$ is source algebra equivalent to $B'_{\mu^w}$, and $B'_{\nu^w}$ is source algebra equivalent to $B_{\mu^w}$.

\end{enumerate}
\end{thm}
\begin{proof}
We have shown in $\S 3$ that if $\nu$ and $\mu$ are extremal points of an $i$-string,       then $K_i$ is a $w$-allowed action, and thus $(\nu,\mu)$ is a $w$-compatible pair.
\begin{enumerate}
\item Suppose that the parities of $\nu$ and $\mu$ are the same.
We have shown that $m_{H,D}({\mathcal O}Gbc)$ is exactly the number calculated in Lemma \ref{K}.  When we sum over $ \Irr(\tilde S_n, b )$ or $ \Irr(\tilde S_m, c )$, which have the same number of elements and get $m_{H,D}({\mathcal O}Gbc)|\Irr(\tilde S_n, b )|$.
Thus,  Theorem 2.5 of \cite{HK} applies and
the block
algebras $ {\mathcal O}\tilde S_nb $  and  $ {\mathcal O}\tilde S_mc $  are source
algebra isomorphic.
\item Suppose the parities are different.  Then this source algebra equivalence is obtained from Lemma 5.1 and Theorem 6.3 of \cite{KS}.
\end{enumerate}
\end{proof}

\section{A sharp bound for Donovan's conjecture}

\begin{defn} We say that two blocks $\rho^{w}$, $\sigma^{w}$ of the same weight are \textsl{allowed equivalent} if one can be obtained from the other by a sequence of $w$-allowed actions.
\end{defn}

Now we wish to find properties which will indicate that a block is allowed-equivalent to a block of lower rank. This will allow us to find a rank $N_0$ such that every block is allowed equivalent to a block of rank $N$, $N\leq N_0$. This was done already in \cite{K}. However, by using crossovers  and by tighter analysis of the possible actions we can make the bound in \cite{K} sharp, and exhibit a block $\rho^w_w$ which attains the bound.

\begin{lem} Let $\rho^{w}$ be a block with the $p$-core $$c(\rho)=((\ell_{1},\epsilon_{1}),...,(\ell_{t},\epsilon_{t}))$$ such that
 for each $i\in I=\left\{1,...,t\right\}$, either $\ell_{i}\geq w$, or $\ell_{i}=0$, then there is a block $\mu^{w}$, of a lower rank, with the $p$-core $c(\mu)=((\ell^{'}_{1},0),...,(\ell^{'}_{r},0),(0,1),...,(0,1))$, that is allowed equivalent, by the $w$-allowed actions, to the block $\rho^{w}$, and such that the
 values of $\ell^{'}_{j}$ form a permutation of those values of $\ell_i$ with $\ell_i \geq w.$
\end{lem}

\begin{proof}
 Let $\rho^{w}$ be a block with the $p$-core \\ $c(\rho)=((\ell_{1},\epsilon_{1}),...,(\ell_{t},\epsilon_{t}))$, satisfying for each $i\in I$ either $\ell_{i}\geq w$, or $\ell_i=0$. (Note that for $\ell_i=0$, by definition, $\epsilon_i=1$).

 \textbf{Step 1:} If all $\epsilon_{i}=0$, then $\rho$ is already in the desired form.  If not, let $k$ be the first place in the $p$-core $\rho$ satisfying $\epsilon_{k}=1$, and let $j$ be the first place in the $p$-core $\rho$, after $k$, satisfying $\epsilon_{j}=0$, if such exists (i.e. all runners from runner $k$ to runner $j-1$ are empty and $\ell_j\neq 0$ i.e. $\ell_j\geq w$).  We are going to run a recursion on $k$ in order to show that we can transform $\rho$ by $w$-allowed actions to the form
 $$c(\nu)=((\tilde{\ell}_{1},0),...,(\tilde{\ell}_{r},0)(\tilde{\ell}_{r+1},1),...,(\tilde{\ell}_{t},1)).$$
 If no $j$ exists, then we can take  $\nu=\rho$ and proceed to the second step.

  If $j$ exists, we have $\ell_{j}+\ell_{j-1}\geq w$ so we can do the $w$-allowed action $K_{j-1}$ and get a block with the $p$-core such that the pair $(\ell_{j-1},0),(\ell_{j},1)$ have been swapped. Currently, the new $\ell_{j-1}$ is the old $\ell_{j}$ and $\ell_{j-1}+\ell_{j-2} \geq w$,  so we can do the $w$-allowed action $K_{j-2}$ and so on. In summary we do the $w$-allowed action $K_{k}\circ K_{k+1} \circ ... \circ K_{j-1}$ and we get to the situation that all runners from runner $k+1$ to runner $j$ are empty.  If there is no $j>k+1$ with $\epsilon_j=0$, we have finished the first step.  Otherwise, we replace $k$ by $k+1$ and continue.

  \textbf{Step 2:}
If, in  the $p$-core
\[c(\nu)=((\tilde{\ell}_{1},0),...,(\tilde{\ell}_{r},0)(\tilde{\ell}_{r+1},1),...,(\tilde{\ell}_{t},1)),\]
all the $\tilde{\ell}_{i}$ for $i>r$ equal $0$, then the lemma has been proven.  If not, let $s$ be the last place in the $p$-core $\rho$ satisfying $\tilde{\ell}_s\neq 0$ meaning $\tilde{\ell}_s\geq w$ and $\epsilon_s=1$.  We do a backwards recursion on $s-r$.

By $K_{t-1}\circ ... \circ K_{s+1}\circ K_s$ we can bring the pair $(\tilde{\ell}_s,1)$ to the place $t$, and perform the $w$-allowed action $K_t$ to invert $\epsilon_t$ from $1$ to $0$. We now repeat Step 1 to get an new $\nu$ with $r$ replaced by $r+1$. When we apply Step 2 again, the new $s'$ will be no greater than the previous $s$, because the actions of Step 1 will return all the pairs which came after $(\tilde{\ell}_s,1)$ to their previous places.

   We get a block $\mu^{w}$  with the $p$-core $$c(\mu)=((l^{'}_{1},0),...,(l^{'}_{r},0),(0,1),...,(0,1))$$ that is allowed equivalent to the block $\rho^{w}$, and of a lower rank than $\rho^{w}$ (because of the $w$-allowed actions that we did  reduce the rank of the $p$-strict partition).  Since we showed in Section 4 that allowed equivalent blocks have equivalent source algebras, and a block is Morita equivalent to its source algebra, we are have actually shown that the blocks are Morita equivalent.

\end{proof}

\begin{lem}\label{rank-lowering} Let $\rho^{w}$ be a block with the $p$-core $$\rho=((\ell_{1},\epsilon_{1}),...,(\ell_{t},\epsilon_{t}))$$ satisfying for each $i\in I=\{1,...,t\}$, either $\ell_{i}> w$ or $\ell_i=0$. There is an allowed equivalent block $\sigma^{w}$ of lower rank with the $p$-core $c(\sigma)=((\ell^{'}_{1}-1,0),...,(\ell^{'}_{r}-1,0),(0,1),...,(0,1))$, where $(\ell^{'}_{1},...,\ell^{'}_{r})$ is a permutation of the non-zero $\ell_{i}$.
\end{lem}

\begin{proof}
Let $(\ell_{i_1},...,\ell_{i_r})$ be the set of all $\ell_i>0$.
First of all, by the previous lemma, the block $\rho^{w}$ is allowed equivalent to a block $\mu^{w}$ with the $p$-core $c(\mu)=((\ell^{'}_{1},0),...,(\ell^{'}_{r},0),(0,1),...,(0,1))$, where $(\ell^{'}_{1},...,\ell^{'}_{r})$ is a permutation of $(\ell_{i_1},...,\ell_{i_r})$.

 We first want to reduce each $\ell_{i}$ ($1 \leq i \leq r$) by one, with $\epsilon_{i}=1$. In terms of the abacus, this will bring all the beads to be on runners $p-r$,...,$p-1$. In order to accomplish this,
 we run a recursion on $i$, for $1\leq i \leq r$.  For $i=1$ we perform $K_0$.
To reduce some $\ell_{i}$ by one we have to bring it to be $\ell_{1}$, by doing $K_{1} \circ K_{2} \circ ... \circ K_{i-1}$ and then do $K_{0}$. We know that $l^{'}_{i}>w$ for $1 \leq i \leq r$ so the involution $K_{0}$ is $w$-allowed action ($l^{'}_{1}-1\geq w$), and also $K_{j}$ when $\epsilon_{j}\neq\epsilon_{j+1}$ is  $w$-allowed action since ($\ell^{'}_{i}+\ell^{'}_{i+1}\geq w$).

We then apply the previous lemma again to change all $\epsilon_{i}$ to $0$.  This is possible because all $\ell_i-1\geq w$.

\end{proof}

\begin{lem} Let $\rho^{w}$ be a block with the $p$-core $$c(\rho)=((\ell_1,\epsilon_1),...,(\ell_t, \epsilon_t)),$$  let $(\ell_{i_1},...,\ell_{i_r})$ be the set of all $\ell_i>0$, let the sequence $(m_{1},...,m_{r})$ be a permutation of  that satisfies $0<m_{1}\leq m_{2}\leq...\leq m_{r}$ and let $m_{gap}$ be the maximal gap between all $m_{i}$.  choose a $j$ such that  $m_{gap}=m_{j+1}-m_{j}$ where $m_{j}>0$ (i.e. $m_{j+1}$ is the smallest of the big $m_{i}$'s and $m_{j}$ is the biggest of the small $m_{i}$'s). If $m_{gap} \geq w$ then all $\ell_{i}$ satisfying $\ell_{i}\geq m_{j+1}$ are reducible by one by $w$-allowed actions.
\end{lem}

\begin{proof}

  We first note that  $m_{j+1} \geq w+m_j>w$. Every pair $(\ell_i,\epsilon_i)$ for which $\ell_i \geq m_{j+1}$, hereafter called a tall pair, can be commuted with every pair $(\ell_k,\epsilon_k)$ for which $\ell_k \leq m_j$, hereafter called a short pair, whether the $\epsilon$ are the same or not, because we always have $\ell_i -\ell_k \leq w$.  Thus if we let $i$ be the first of the tall pairs when the runners are ordered from $1$ to $p-1$, we can move it toward  the front by $w$-admissible actions of type $K_{i'}$ for $0<i'<t$. If $\epsilon_i =1$, then when we reach runner $t+1$, we must perform the action $K_t$, but this is also admissible since
$2\ell_i+1 \geq w$.  Performing these actions recursively, we reach a situation in which all the tall pairs have $\epsilon = 0$ and are in $1$ through $s$, for some $s<r$.
Then they can all be lowered by one as in Lemma \ref{rank-lowering}.

\end{proof}

\begin{prop}\label{maxrank}  Let $\rho^{w}$ with $c(\rho)=((\ell_{1},\epsilon_{1}),...,(\ell_{t},\epsilon_{t}))$ be a block that cannot be reduced by $w$-allowed actions, and let $m_{gap}$ be as before. Then
\begin{itemize}
\item the $p$-core  $ \rho=$ satisfies $min\{\ell_i|1\leq i\leq r\}\leq w$.
\item $m_{gap}\leq w-1$ i.e. the maximal gap is $w-1$.  Among these blocks, that with maximal rank has core $t$-tuple of the form
 $ c(\rho_w)=((w,0),...,(w+(w-1)(t-1),0))$.
\end{itemize}
\end{prop}

\begin{proof}
By Lemma 5.1 we get that every core satisfying, for each $i\in I=\{1,...,t\}$ either $\ell_i\geq w$ or $\ell_i=0$  can be reduced by $w$-allowed actions to a core of the form
$((\ell_1,0),...,(\ell_r,0),(0,1),...,(0,1))$.

By Lemma 5.2 a core $((\ell_1,0),...,(\ell_r,0),(0,1),...,(0,1))$ satisfying \\ $\ell_i>w$ for $1\leq i \leq r$ can be reduced to a core \\ $((l^{'}_1-1,0),...,(l^{'}_r-1,0),(0,1),...,(0,1))$ i.e. in the core that cannot be reduced there is $i$ satisfying $0\neq \ell_i\leq w$.

By Lemma 5.3 if there is a gap, $m_{gap}=m_{j+1}-m_j$ satisfying $m_{gap}\geq w$ then this gap can be reduced until it less then $w$. i.e. in the core that cannot be reduced, the maximal gap is $w-1$.
The maximal rank is attained when the minimum is as large as possible, all gaps are as large as possible, and the ordering of the runners give the largest possible rank.  This give the core $\rho_w$ of statement of the lemma.
\end{proof}

\begin{remark} A block similar to our $\rho_w^w$ was used by Chuang and Kessar \cite{CK} to complete the Chuang-Rouquier proof of the Brou\'e conjecture for the symmetric groups \cite{CR}. It is sometimes referred to in the literature as a RoCK-block, and we will call $\rho_w$ a RoCK-core.
\end{remark}

\begin{lem} Let $c(\nu)=((\ell_1,\epsilon_1),...(\ell_t,\epsilon_t))$ be a core. The rank of this block is
\[N(\nu)=\sum_{i=1}^{t} \ell_i\cdot i^{1-\epsilon_i}\cdot(p-i)^{\epsilon_i}+\frac{\ell_i(\ell_i-1)}{2}\cdot p.\]
\end{lem}

\begin{proof}

 For every pair $(\ell_i,\epsilon_i)$ we consider the addition made to the rank by all the beads on the $i$-th runner. If $\epsilon_i=0$ then there are $\ell_i$ beads, corresponding to parts of the form $ap+i$ for $a=0,...,\ell_i-1$, and if $\epsilon_i=1$ then there are $\ell_i$ beads, corresponding to parts of the form $ap+p-i$ for $a=0,...,\ell_i-1$.

 If $\epsilon_i=0$ then the rank will be
\[\ell_i \cdot i + \frac{\ell_i(\ell_i-1)}{2}\cdot p\] and if $\epsilon_i=1$ then that rank will be
\[\ell_i \cdot (p-i) + \frac{\ell_i(\ell_i-1)}{2}\cdot p.\] So in each case we get that the rank is
\[ \ell_i\cdot i^{1-\epsilon_i}\cdot (p-i)^{\epsilon_i}+\frac{\ell_i(\ell_i-1)}{2} \cdot p \]
for  every pair $(\ell_i,\epsilon_i)$, and we finished.
\end{proof}

\begin{thm}\label{bound} The block $\rho^w_w$ of the maximal rank $N$ which does not lie at the maximal rank  end of any $i$-string has rank
\[N=pw+(\frac{p(w-1)}{2}+1)\cdot(\sum^{t}_{i=1}i^2(w-1)+i)\]
\end{thm}

\begin{proof}

 The block of maximal rank according to \ref{maxrank} before has to fulfill the following:
\[\forall i, \epsilon_i=0\]
\[\ell_1=w\]
\[\ell_2=2w-1\]
\[\ell_3=3w-2\]
\[\vdots\]
\[\ell_t=tw-(t-1)\]
Now we substitute these values in the formula of the previous lemma:

\[ N(\rho_w)=\sum^{t}_{i=1}(iw-(i-1))i+\frac{(iw-(i-1))(iw-i)p}{2},  \]

so

\[ N(\rho_w)=\sum^{t}_{i=1}(i^2w-(i-1)i)+\frac{p(w-1)}{2}\sum^{t}_{i=1}(i^2w-(i-1)i), \]

and we get

\[ N(\rho_w)=(\frac{p(w-1)}{2}+1)\cdot(\sum^{t}_{i=1}i^2(w-1)+i).  \]
Finally, we add $pw$ for the weight of the block:
\[ N(\rho_w^w)=pw+(\frac{p(w-1)}{2}+1)\cdot(\sum^{t}_{i=1}i^2(w-1)+i).  \]
\end{proof}

\begin{remark}
The number calculated in \ref{bound} is actually an integer.
\begin{itemize}
\item If w is odd, then $w-1$ is even, so $\frac{w-1}{2}$ is an integer.
\item If w is odd, then the parity of the sum depends on the parity of
$\sum^{t}_{i=1} i-i^2$, and each term of this sum is even.
\end{itemize}
\end{remark}

\section{Representations of affine Lie algebras}

The bound given in Section 5 was derived from a study of the block-reduced crystal graph, which, in turn, is related to the representations of the twisted affine Lie algebra $A^{(2)}_2t$.  The untwisted affine Lie algebra $A^{(1)}_{\ell}$ has a  diagram which is the Dynkin diagram of the classical Lie algebra $A_\ell$ with an added point labeled $0$ joining the two ends.  The twisted affine algebra is the fixed algebra under the graph automorphism sending $1$ to $\ell$, $2$ to $\ell-1$, etc., leaving $0$ fixed.  We will be using this twisted algebra in the case where $\ell$ is $p-1$ for an odd prime $p$.

    Every affine Lie algebra has a diagram which is a one-point extension of a classical Lie algebra.  The classical Lie algebra for $A^{(1)}_{\ell}$ is of course $A_\ell$.
The classical Lie algebra for $A^{(2)}_2t$ is of type $C_t$. The Cartan matrix $C$ for the twisted affine algebra with $p=11$ is given below, and the Cartan matrix $C_0$ for the corresponding classical algebra is obtained by crossing out the first row and column.

The extended Dynkin diagram is given by
$$ \cdot \Leftarrow \cdot - \cdot - \cdot - \cdot - \cdot \Leftarrow \cdot$$

$$C=
\begin{bmatrix}
2&-2&0&0&0&0\\
-1&2&-1&0&0&0\\
0&-1&2&-1&0&0\\
0&0&-1&2&-1&0\\
0&0&0&-1&2&-2\\
0&0&0&0&-1&2
\end{bmatrix}
$$

The Lie algebra $\mathfrak{G}$ is determined by a Chevalley basis $\{e_i, f_i, h_i| i = 0, \dots,t\}$ over $\mathbb{C}$. All elements of the basis are eigenvectors for the action of the abelian Cartan subalgebra $\mathfrak H$ generated by the $\{h_i\}$, and the elements $a_{ij}$ of the Cartan matrix give the eigenvalues for the $h_i$ acting on $e_j$, with the eigenvalues for the $f_j$ being the negatives.
    The set of \textsl{weights} of the algebra are the possible eigenvalues associated with the $h_i$ for all possible one-dimensional representations.  The weights of the various one-dimensional representations in the adjoint representation of $\mathfrak H$ acting on the entire Lie algebra are called the \textsl{roots}.   The columns $\a_j$ of the Cartan matrix provide a spanning set for the set of roots and are called the \textsl{simple roots}. Every other root can be obtained as a sum of simple roots (giving a \textsl{positive root}) or a sum of their negatives (giving a \textsl{negative root}). The elements of the basis of the Cartan subalgebra are called the \textsl{simple coroots}.   The \textsl{weight space} is the real vector space $V$ generated by the dual basis $\Lambda_0,\dots,\Lambda_\ell$ to the $h_0,\dots,h_\ell$. The columns $\a_j$ of the Cartan matrix generate an integral lattice in $V$ called the root lattice.

    The determinant of the Cartan matrix of the affine Lie algebra is zero, the rank being one less than the degree.  The linear combination of roots which is in the kernel of the action of $\mathcal H$ on the Lie algebra is called the \textsl{imaginary root} and is generally denoted by $\delta$. In our case it is given by:
    $$\delta = \a_t+2\a_{t-1}+\dots+2\a_0.$$
    There is a corresponding linear combination of coroots which is central in the Lie algebra, i.e., has a trivial bracket with every element.  In our case it is given by
$$c= 2\a_t+\dots+2\a_1+\a_0.$$

    The Cartan matrices of the affine Lie algebras can be \textsl{symmetricized}, by multiplying the basis elements by rational constants.  The symmetricized matrix $B$ can then be used to determine a bilinear form.

$$B=
\begin{bmatrix}
1&-1&0&0&0&0\\
-1&2&-1&0&0&0\\
0&-1&2&-1&0&0\\
0&0&-1&2&-1&0\\
0&0&0&-1&2&-2\\
0&0&0&0&-2&4
\end{bmatrix}
$$

      In the symmetrized form $(,)$ for our twisted affine algebra of type $A$, the lengths of the simple roots are
    $1,2, \dots,2,4$.  Of particular importance are the ``long" roots of the classical algebra, which in our case consist of all roots of length $4$. We define
    $$\beta_i = \a_t + 2\a_{t-1} + \dots + 2\a_i, i=1,\dots,t.$$ Then
    $$(\beta_i,\beta_i) = (\a_t,\a_t)+[2(2\a_{t-1},\a_{t})+(2\a_{t-1},2\a_{t-1}]$$ $$+\sum^{t-1}_{j=i+1} [2(2\a_i,2\a_{i+1})+(2\a_i,2\a_i)]$$
    $$=4-8+8+\sum^{t-1}_{j=i} ((-8)+8)=4.$$
These are all the long roots.

    The importance of these long roots is in determining the infinite Weyl group. The Weyl group is the group of automorphisms of the root lattice, and is generated by reflections in the simple roots, the reflection $r_\a$ in a root $\a$ being the involution fixing the hyperplane perpendicular to $\a$ and reversing the coordinate in the direction of $\a$, so that $\a$ goes to $-\a$ and more generally,
    $r_\a(\beta) = \beta - 2\frac{(\beta,\a)}{(\a,\a)}\a$.

    The Weyl group of an affine Lie algebra is a semidirect product with a quotient isomorphic to the finite Weyl group $W_0$ of the classical Lie algebra $\mathfrak G_0$ and an abelian normal subgroup $T$ which has $\ell$ generators corresponding to the long roots of the classical algebra.  The elements of $T$ are labeled by weights and satisfy an equation $t_\a \circ t_\gamma= t_{\a+\gamma}$, which demonstrates their commutativity. The formula for $t_\a$ acting on the weight $\Lambda_0$ is
    $$t_\a(\Lambda_0) = \Lambda_0 + \a - \frac{1}{2}|\a|^2\delta.$$  In our case, the generators are the elements $t_{\frac{1}{2}\beta_i}$, for $i=1,\dots,t$.

    There is a second way to understand the Weyl group, given in Exercise 6.7 of \cite{Ka}.  The group is represented by actions $s_0, s_1, \dots, s_t$ on the real vector space of coroots of $\mathfrak G_0$, generated by $h_1, h_2, \dots h_t$.  Each $s_i$ for $i>0$ acts as a reflection in the plane defined by $\alpha_i = 0$, (that is, the set of all coroots whose product with $\alpha_i$ in the symmetric bilinear form is 0). The generator $s_0$ corresponds to a reflection in the plane $\beta_1 = \frac{1}{2}$.

    There is a representation of our Lie algebra $A^{(2)}_2t$ on the complexified sum of the Grothendieck groups of the group algebras of $\tilde S_n$ and $\tilde A_n$, for all n, with copies of K for $n=0$. The representation is given by sending $e_i$ to a restriction operator, which first restricts the simple to a group of rank one lower and then to the correct block.  The top can then be shown to correspond exactly to the simple determined by removing the $i$-good node.  However, this must be done, as usual, with crossovers whenever the parities differ, as they always do when $i \neq 0$.  This representation of $A^{(2)}_2t$ corresponds to the highest weight representation generated by the fundamental dominant weight $\Lambda_0$, the dual basis vector corresponding to $h_0$.  The weights of all the simples corresponding to a block with content $\gamma= (\gamma_0,\dots,\gamma_t)$ are given by
    $$\Lambda_0 - \sum^t_{i=0}\gamma_i\a_i.$$  By \cite{Ka}, this set of weights are all of the form $$\{w\cdot\Lambda_0 -k\delta|w \in W, k = 0,1,\dots\}.$$  If we let
    $$X_k=\{w\cdot\Lambda_0 -k\delta|w \in W\}, k = 0,1,\dots,$$ then the set $X_0$, called the set of maximal weights, is the set of weights  from which no copies of $\delta$ can be removed. These correspond to the blocks of defect $0$.  More generally, $k$ corresponds to the weight of the block, and the elements $\{w \cdot \Lambda_0|w \in W\}$, to the core.  Now, in fact, $\Lambda_0$ is fixed under all the elements of $W_0$, which operates only on elements of the classical Lie algebra.  Therefore, the set of weights corresponding to cores is of the form $T \cdot \Lambda_0$.  Since, as we have already noted, the elements of $T$ correspond to elements of the lattice
    $$\frac{1}{2}(n_1\beta_1+\dots n_t\beta_t),$$ we need to understand the transformation between the integer t-tuple $(n_1,\dots,n_t)$ and the core $t$-tuple $((\ell_1,\epsilon_1),\dots,(\ell_t,\epsilon_t))$.

\begin{lem} The element with coordinate vector $\bar n = (n_1,\dots,n_t)$ in the half-integer weight lattice of the long roots $\{\beta_1,\dots,\beta_t\}$ corresponds to the core satisfying $$n_i = (-1)^{\epsilon_i}\ell_i.$$
\end{lem}
\begin{proof}
     We begin by calculating the core corresponding to $\frac{1}{2}n_1\beta_1$. Since $\beta_1 = \delta -2\a_0$,  we have a core $\rho$ whose weight is
    $$t_{\frac{1}{2}n_1\beta_1}(\Lambda_0)= \Lambda_0 + \frac{1}{2}n_1\beta_1 - \frac{1}{2}n_1^2\delta$$
    $$=\Lambda_0 -(n_1^2-n_1)\delta-n_1\a_0.$$
    For $n_1>0$ this is precisely the weight of the core with $n_1$ beads on runner $1$, and for $n_1<0$, it is the weight of the core with $|n_1|$ bead on runner $p-1$, The two cases being gotten by adding or removing zeros from the Young diagram of a stack of beads $n_1$ on the $0$-runner, which contributes $1+2+\dots+n_1$ copies of $-\delta$ to the weight.

    We now operate on the half-integer lattice by simple reflections $r_{\a_i}$, for $0<i<t-1$.  We have
    $$r_{\a_i}(\beta_i) =\beta_i-2\frac{(\beta_i,\a_i)}{2}\a_i$$
    $$= \beta_{i+1} + 2\a_i - (2(\a_{i+1},\a_i)+2(\a_i,\a_i))\a_i$$
    $$=\beta_{i+1} +2\a_i -(-2+4)\a_i$$
    $$=\beta_{i+1}.$$
    Since $r_{\a_i}$ is an involution, we see that it must carry $\beta_{i+1}$ to $\beta_i$. The calculation in the special case $i=t-1$ is slightly different, because there is only one copy of $\a_t$ in $\beta_{t-1}$, but leads to the same result because $\a_t$ is a long root.
    We can then apply the reflections $r_{\a_i}$ for $i=1,2,\dots,t-1$, to move the stack of beads along the runners.  To finish the proof, we use the fact that we have a homomorphism of $\mathbb Z^t$ onto $T$.

\end{proof}
\begin{corl}
     The actions $K_i$ on the cores correspond to reflections $r_{\a_i}$ in the simple roots.
\end{corl}
\begin{proof}
    \begin{enumerate}
    \item For $0 <  i < t$, the action $K_i$ switches the pairs $(\ell_i,\epsilon_i)$ and and  $(\ell_{i+1},\epsilon_{i+1})$, while $r_{\a_i}$ switches $\beta_i$ with $\beta_{i+1}$, thus switching $n_i$ and $n_{i+1}$.
    \item For $i=t$, the action of $K_t$ sends $\epsilon_t$ to $1-\epsilon_t$, while $r_{\a_t}$ sends $\beta_t$ to $-\beta_t$.
    \item  For $K_0$, both the parity changes and the number of beads is raised or lowered.  On the Lie algebra side, the situation in more complicated than before, because $\a_0$ does not lie in the hyperplane containing the $\beta_i$. After acting by $\a_0$, one must eliminate the copies of $\a_0$ which appear by subtracting off an appropriate multiple of $\delta$.  The total result is to send $n_1\beta_1$ to $\beta_1-n_1\beta_i$, which is a reflection around the midpoint of $\beta_1$, as described after the definition of the Weyl group above.  This is the operation which allows the set of weights to become infinite, as for example, $\bar n = (-1,0,\dots,0)$ is reflected to $(2,0,\dots,0)$.
\end{enumerate}
\end{proof}

We will use this information to bound the number of possible Morita equivalence classes of blocks of a given weight.  Let us define the \textsl{level} of a block $\rho^w$ to be the value of $\gamma_t(\rho^w)$. In terms of the coordinates $(n_1,\dots,n_t)$, the formula is given by $$\sum_{i=1}^t \frac{n_i^2-n_i}{2}.$$  This level number is invariant under all permutations of the $n_i$, and also under all maps sending $n_i$ to $1-n_i$.
By the work done in \cite{AS}, we know that all the cores satisfying the conditions in Proposition \ref{bound} lie in levels less than or equal to that of the ''RoCK-block" $\rho_w^w$,with the  core coordinates $n_i = w+(i-1)(w-1)$.  Since each rise of $1$ in the exponent adds a copy of $\delta$ to the content, and the $t$-coordinate of $\delta$ is one, the exponent of $w$ just shifts the level up by $w$.  Thus the number of cores which are candidates to be representatives of distinct Morita equivalence classes of blocks of weight $w$ is the number of cores $\rho$ with level $$\gamma_t(\rho) \leq \gamma_t(\rho_w).$$  For $w>0$, this level number, by the formulae in Proposition \ref{bound}, is
$$L_w=\sum_{i=1}^t (w +(w-1)(i-1))((w-1)+(w-1)(i-1))$$
$$=\sum_{i=1}^t (w +(w-1)(i-1))(w-1)i.$$  The set of all $\bar n$ with level less than or equal to a level $L$ is roughly shaped like a $t$-sphere centered on $\frac{1}{2}(1,1,\dots,1)$ with radius of the order of magnitude $\sqrt{L}$.

  The number of interior integer points is of the order of magnitude $L^\frac{t}{2}$.  In particular, since $L_w$ is of the order of $t^2w^2$,  the number of possible candidates for Morita equivalence class representatives is of the order of $t^tw^t$, and thus, for a fixed prime $p$, is polynomial in $w$ of degree $t$.
\begin{examp}
In the matrix below, we give the levels for $p=5$ for blocks with $n_i$ in the range $-4 \leq n_i \leq 5$.  There are four axes of symmetry, all passing through the point $( \frac{1}{2},\frac{1}{2})$ at the center of the matrix, generated by reflections in the simple roots $\a_0$ and $\a_1$.  The reflections in the simple root $\a_2$ are the reflections whose axis is the lower of the two rows containing $0$.  These reflections change the level by a fixed amount, determined by the row being reflected, and corresponding to the coordinate $n_2$, which counts the number of beads on runner $2$ or $3$.
$$
\begin{bmatrix}
20&16&13&11&10&10&11&13&16&20\\
16&12&9&7&6&6&7&9&12&16\\
13&9&6&4&3&3&4&6&9&12\\
11&7&4&2&1&1&2&4&7&11\\
10&6&3&1&0&0&1&3&6&10\\
10&6&3&1&0&0&1&3&6&10\\
11&7&4&2&1&1&2&4&7&11\\
13&9&6&4&3&3&4&6&9&12\\
16&12&9&7&6&6&7&9&12&16\\
20&16&13&11&10&10&11&13&16&20

\end{bmatrix}
$$
\end{examp}

\end{document}